\documentclass[%
 notitlepage,
 amsmath,amssymb,
 aps,
]{revtex4-2}

\usepackage{graphicx}
\usepackage{dcolumn}
\usepackage{bm}
\usepackage{hyperref}
\usepackage{placeins}
\usepackage{comment}

\begin{document}

\preprint{APS/123-QED}

\title{A $p$-adaptive high-order mesh-free framework for fluid simulations in complex geometries}

\author{Ruofeng Feng}
 \email{sdufrf@163.com}
\affiliation{%
Department of Mechanical and Aerospace Engineering, University of Manchester, UK
}%
\author{Jack R. C. King}%
\email{jack.king@manchester.ac.uk}
\affiliation{%
Department of Mechanical and Aerospace Engineering, University of Manchester, UK
}%
\author{Steven J. Lind}
\affiliation{%
School of Engineering, Cardiff University, Cardiff, UK
}%

\date{\today}

\begin{abstract}
This paper presents a $p$-adaptive, high-order mesh-free framework for the accurate and efficient simulation of fluid flows in complex geometries. High-order differential operators are constructed locally for arbitrary node distributions using linear combinations of anisotropic basis functions, with weights constructed to reproduce polynomials up to a prescribed $p$ order. A dynamic $p$-refinement strategy is developed to refine (increase) or de-refine (decrease) the polynomial order used to approximate derivatives at each node. A new refinement indicator for mesh-free methods is proposed, based on local error estimates of the Laplacian operator, and is incorporated into the solution procedure at minimal added computational cost. Based on this error indicator, a refinement criterion is established to locally adjust the polynomial order $p$ for the solution. The proposed adaptive mesh-free scheme is then applied to a range of canonical PDEs, and its potential is demonstrated in two-dimensional simulations of a compressible reacting flow in porous media. For the test cases studied, the proposed method exhibits potential to save up to 50\% of computational costs while maintaining the specified level of accuracy. The results confirm that the developed $p$-adaptive high-order mesh-free method effectively captures highly non-linear regions where high-order approximation is necessary and reduces computational costs compared to the non-adaptive method, preserving high accuracy and solution stability.

\begin{description}
\item[Keywords]
High-order; Mesh-free; Error-based indicator; $p$-adaptivity; Fluid simulation; Porous media
\end{description}
\end{abstract}

\maketitle


\section{\label{sec:level1}Introduction}

Simulating fluid problems via the numerical solution of partial differential equations (PDEs) relies on the accurate calculation of spatial derivatives from field values at a discrete set of points. High-order approximations are particularly desirable, as they can achieve extremely high accuracy with lower resolution, thereby reducing computational costs. High-order finite difference (FD)\cite{lele1992} and spectral methods\cite{karniadakis1989}, for instance, are easily constructed on structured meshes but they are consequently limited to simple geometries. While immersed boundary methods can extend these approaches to complex geometries, they face significant challenges in accurately imposing boundary conditions\cite{rauch2018}. Conversely, methods such as unstructured finite volume (FV)\cite{weller1998} and finite element methods address the challenge of handling complex geometries. Their primary difficulty, however, is that accuracy is highly dependent on mesh quality. In such cases, generating a high-quality, body-fitted unstructured mesh with minimal skewed elements becomes a major bottleneck. This process can be so resource-intensive, in some cases, such as porous media, taking longer than the simulation itself, that it undermines the efficiency of the method\cite{wood_2020}. 

Mesh-free methods offer an alternative to traditional grid-based approaches by operating on a set of computational points (often called particles) without requiring information on topological connectivity. A key advantage of mesh-free methods is that discretising a complex geometry with an unstructured set of nodes is straightforward and easily automated (e.g.,\cite{fornberg_2015a}). Many of these methods, particularly strong-form variants, work by constructing differential operators on a local stencil of collocation points. They have several variants, including radial basis function (RBF) methods~\cite{kansa_1990a,kansa_1990b,shu_2003,cecil_2004,wright_2003,fornberg_2015}, generalised finite difference methods (GFDM)\cite{jensen_1972,benito_2007,gavete_2017}, generalised moving least squares (GMLS) methods~\cite{trask_2017,trask_2018}, reproducing kernel particle methods (RKPM)~\cite{liu_1995,liu_1995b}, and finite particle methods (FPM)~\cite{zhang_2004,liu_2005,asprone_2010,asprone_2011}. The specific order of approximation can be constructed using various basis functions (e.g., orthogonal polynomials, Taylor monomials) and techniques (e.g., weighted or moving least squares). A primary challenge for these collocation-based methods is ensuring the conservation of physical quantities (such as mass and momentum), a difficulty inherent to the strong-form formulation. In this regard, Smoothed Particle Hydrodynamics (SPH) stands out as a popular alternative that inherently preserves mass (and momentum depending on the formulation) (\cite{lind_2020}) and has been used in a wide range of applications (\cite{feng2021,dominguez2019,king2023,khayyer2024,feng2024,schaller2024}). However, a major limitation of SPH is its accuracy, in which traditional Lagrangian SPH is often restricted to low-order convergence rates as resolution increases. In recent years, high-order formulations of SPH have been widely explored to improve the consistency, such as  consistency corrections to kernels~\cite{lind2016} and derivative operators~\cite{nasar2021}, Riemann-type formulations\cite{vila1999,parshikov2002}, W/TENO reconstruction\cite{meng2024}, and a hybrid spectral scheme\cite{lin2025}. The Local Anisotropic Basis Function Method (LABFM), proposed by King et al.\cite{king2020high}, is another variant of high-order mesh-free methods. It constructs derivative operators from a weighted sum of local field value differences, where the weights are formed from a linear combination of Anisotropic Basis Functions (ABFs) to guarantee polynomial consistency of a specified order. LABFM may be formulated at higher orders (commonly up to 10th order) and has a lower computational cost for a given level of accuracy compared to other consistent mesh-free methods\cite{king2022high}. Nevertheless, for the practical application of high-order mesh-free methods, the major obstacle remains their intensive computational cost due to large stencils to achieve higher-order consistency.

Adaptive solution strategies are an important approach for balancing computational efficiency and high accuracy. These techniques are particularly effective for problems where high resolutions are required only in localised regions, such as simulating the dynamic evolution of fine flame fronts. Two primary adaptive approaches have been proposed in the literature: $p$-adaptivity, and $h$-adaptivity. In $p$-adaptivity, the accuracy is varied by locally adjusting the order of the interpolant approximation. In contrast, $h$-adaptivity locally adjusts the resolution of the spatial discretisation. $h$-adaptivity has seen widespread adoption in mesh-based methods for fluid dynamics, commonly known as adaptive mesh refinement (AMR) techniques. Notable examples include Basilisk\cite{popinet2015}, which utilises a quadtree/octree-based grid structure with a second-order accurate FV discretisation, the AMReX library\cite{zhang2019amrex}, which employs a block-structured approach to manage nested patches of successively finer grids, and the HAMISH code\cite{cant2022}, which combines adaptive Cartesian grids with a high-order scheme using fourth-order flux reconstruction methods. $p$-adaptivity, on the other hand, has been frequently explored within high-order numerical frameworks such as spectral element methods (SEM)\cite{proot2002,zayernouri2014} and discontinuous Galerkin (DG) methods\cite{naddei2019comparison,mossier2022} to improve the convergence rate. $hp$-adaptivity is another adaptive strategy, e.g.,\cite{mitchell2011,basile2022}, that integrates $h$- and $p$- refinement. With the incorporation of the $hp$-decision strategy often based on the measure of solution smoothness, the method allocates $p$-refinement to regions with a smooth solution, while $h$-refinement is applied to non-smooth regions.

In the context of mesh-free methods, $h$-adaptivity can be implemented by adding or removing computational nodes \cite{slak2019adaptive}, without the issues associated with mesh-based methods (e.g., topological data update \cite{ji2010}, hanging nodes\cite{lian2006}, and mesh quality\cite{gou2017}). $h$-adaptivity has been studied extensively for various mesh-free methods, such as GFDM\cite{benito2003h}, RBF-FD\cite{oanh2017adaptive}, GMLS\cite{hu2019spatially}, and SPH\cite{vacondio2013variable}. While considerable research has focused on $h$-adaptive mesh-free methods, few studies have been devoted to $p$-adaptivity. The very recent study by Jančič and Kosec \cite{janvcivc2024} has achieved $p$-refinement in RBF-FD for the linear elasticity problem. However, due to the need for explicit re-evaluations of the PDEs at a higher order approximation for the refinement indicator, the adaptive solution may become more expensive than the non-adaptive method.

In light of these challenges, this work presents a $p$-adaptive high-order mesh-free scheme built upon LABFM. A novel refinement indicator is proposed, which is derived from a local estimate of the error associated with the mesh-free approximation of second-order differential operators. This indicator is then used to establish a procedure that locally adapts the polynomial order, $p$, according to user-defined error thresholds, with low additional computational costs.  The performance and accuracy of the proposed $p$-adaptive LABFM are evaluated on a range of benchmark tests, including a standard test function convergence study, the viscous Burgers’ equation, and the Kelvin-Helmholtz instability. Finally, the proposed method is applied to simulate compressible reacting flows, investigating the propagation of hydrogen–air flames in both unconfined and confined geometries. The computational speed-up achieved through $p$-adaptivity is discussed. To the best of the authors’ knowledge, this is the first time that $p$-adaptivity has been enabled within a mesh-free framework for fluid simulations.

The remaining part of this paper is organised as follows: Section~\ref{sec:labfm} presents a brief description of the LABFM methodology for high-order spatial discretisation. In Section~\ref{sec:p-adapt}, the formulation of the refinement indicator as well as the $p$-adaptivity strategy is introduced. The accuracy and performance of the developed method are examined in Section~\ref{sec:tests} through a number of numerical benchmark tests. Finally, conclusions are summarised in Section~\ref{sec:conclusions}.

\section{High-order mesh-free scheme}
\label{sec:labfm}
The high-order mesh-free scheme adopted in this study follows the Local Anisotropic Basis Function Method (LABFM)~\cite{king2020high}. LABFM has demonstrated up to 10th-order convergence, even for node distributions with significant local disorder \cite{king2022high}. This allows the method to remain extremely accurate, even in complex geometries where non-uniformity in the node distribution is unavoidable. 

\subsection{Evaluation of spatial derivatives}
To start with, we define a general expression for the spatial derivative as follows,
\begin{equation}
\label{eq:Ld0}
L^{d}(\cdot) = \mathbf{C}^{d} \cdot \mathbf{D}(\cdot)
\end{equation}
where \( d \) is the index that identifies the partial derivative under consideration, \(\mathbf{D}(\cdot)\) is the vector of partial derivatives that is expressed as,
\begin{equation}
\label{eq:D()}
\mathbf{D}(\cdot)\ = \begin{bmatrix} 
\frac{\partial(\cdot)}{\partial x}, & \frac{\partial(\cdot)}{\partial y}, & \frac{\partial^2(\cdot)}{\partial x^2}, & \frac{\partial^2(\cdot)}{\partial x\partial y}, & \frac{\partial^2(\cdot)}{\partial y^2}, &  
\frac{\partial^3(\cdot)}{\partial x^3}, & \frac{\partial^3(\cdot)}{\partial x^2\partial y}, & \frac{\partial^3(\cdot)}{\partial x\partial y^2}, & \frac{\partial^3(\cdot)}{\partial y^3}, & \frac{\partial^4(\cdot)}{\partial x^4}, 
\cdots
\end{bmatrix}^T
\end{equation}
in which \(\mathbf{C}^{d}\) is a coefficient vector selecting the required derivatives for a given index \( d \). For example, to estimate the gradients and Laplacian of a function, for which case \(L^{x} = \partial/{\partial x}\), \(L^{y}=\partial/{\partial y}\), \(L^{L}=\nabla^{2}\), \(\mathbf{C}^{d}\) can be expressed as,

\begin{equation}
\mathbf{C}^d = 
\begin{cases} 
\begin{bmatrix} 
1, 0, 0, 0, 0, 0, \dots 
\end{bmatrix}^{T} & \text{if } d = x \\[2ex]
\begin{bmatrix} 
0, 1, 0, 0, 0, 0, \dots  
\end{bmatrix}^{T} & \text{if } d = y \\[2ex]
\begin{bmatrix} 
0, 0, 1, 0, 1, 0, \dots  
\end{bmatrix}^{T} & \text{if } d = L 
\end{cases}
\end{equation}

The mesh-free construction of spatial derivatives in Eq.~(\ref{eq:Ld0}) to polynomial consistency of order $p$ follows the methodology of LABFM. In the following, only the essential components of LABFM are presented. For a detailed description of the methodology, readers are referred to \cite{king2020high,king2024cmame}.

\begin{figure}[h]
    \centering
    \includegraphics[width=0.49\textwidth]{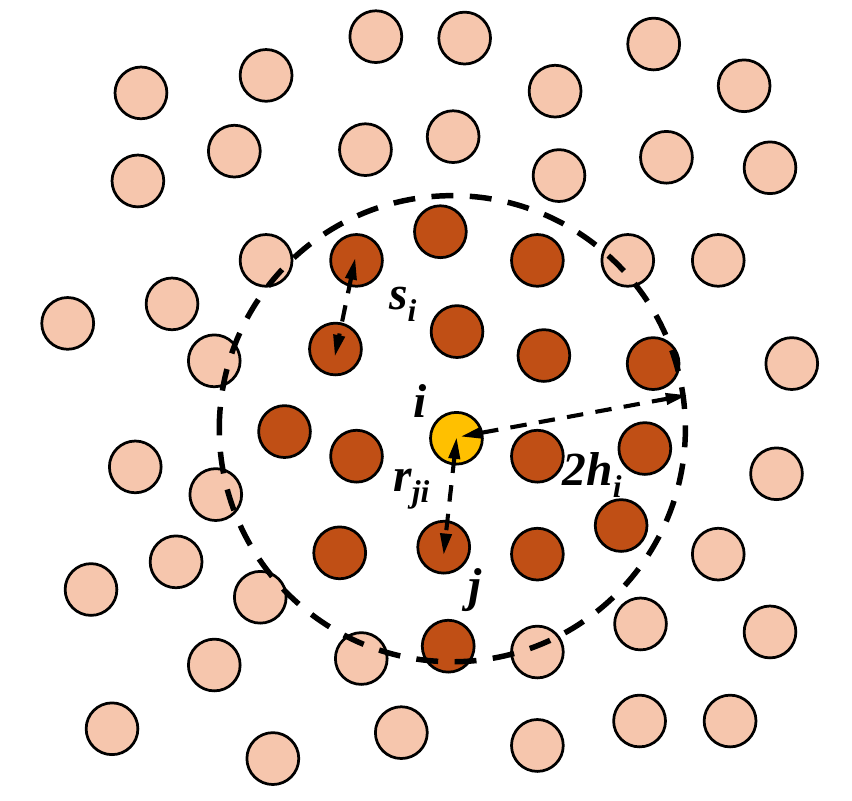}
    \caption{An illustration of the computational stencil.} \label{fig:stencil}
\end{figure}

In LABFM, the computational domain is discretised as a set of $N$ unstructured nodes, as illustrated in Fig.~\ref{fig:stencil}. Each node $i \in [1,N]$ has position vector $\bm{r}_{i} = (x_{i}, y_{i})^{T}$, a distribution length scale $s_{i}$, and a characteristic stencil length scale $h_{i}$. The distribution length scale $s_{i}$ represents the local average inter-node distance and serves as a resolution metric. The resolution may vary spatially, allowing for non-uniform spatial adaptivity. The computational stencil for node $i$ contains $\mathcal{N}_{i}$ neighbouring nodes $j$ satisfying $|\bm{r}_{ji}| \leq 2h_{i}$. For a two-dimensional stencil with  uniform node density, the neighbour count approximately follows $\mathcal{N}_{i} \approx 4\pi(h_{i}/s_{i})^{2}$, with the ratio $h_{i}/s_{i}$ characterising stencil size.  Property differences between nodes are defined as $(\cdot)_{ji} = (\cdot)_{j} - (\cdot)_{i}$. 

The general discrete operator at a computational node \(i\) to reproduce the spatial derivative in Eq.~(\ref{eq:Ld0}) is given by,
\begin{equation}
\label{eq:Ld}
L^{d}(\cdot)_{i,h,p} = \sum_{j} (\cdot)_{ji} w_{ji}^{d}
\end{equation}
where the subscripts \(h\) and \(p\) represent the local stencil length scale and local polynomial order adopted in the interpolation, respectively. \( w_{ji}^d \) are the set of weights. Here (and throughout), the summation over \( j \) is over all $\mathcal{N}_{i}$ nodes within the stencil of node $i$. 

In Eq.~\eqref{eq:Ld}, the weights \( w_{ji}^d \) in LABFM are constructed as the weighted sum of a series of Anisotropic Basis Functions (ABFs) \(\mathbf{W}_{ji} = \mathbf{W}(\bm{r}_{ji} / h_{i}) \), given by
\begin{equation}
w_{ji}^{d} = \mathbf{W}_{ji} \cdot \mathbf{\Psi}^{d}{}_{i}
\end{equation}

Following~\cite{king2022high} the ABFs are constructed from bi-variate Hermite polynomials multiplied by a radial basis function (RBF) \(\psi\), with the $\alpha$-th component of $\mathbf{W}_{ji}$ expressed as,
\begin{equation}
W_{ji}^{\alpha} = \frac{\psi(r_{ji}/h_{i})}{\sqrt{2^{\alpha_1+\alpha_2}}} \, H_{\alpha_1}\left(\frac{x_{ji}}{h_{i}\sqrt{2}}\right) \, H_{\alpha_2}\left(\frac{y_{ji}}{h_{i}\sqrt{2}}\right),
\label{eq:abf}
\end{equation}
where $r_{ji}=|\bm{r}_{ji}|$ is the distance between nodes $i$ and $j$, $H_{a}$ is the $a$-th order univariate Hermite polynomial, the radial basis function \(\psi\) adopted in this work is the Wendland C2 kernel \cite{dehnen2012improving}, $\alpha_{1}$ and $\alpha_{2}$ are the $x$ and $y$ exponents of the $\alpha$-th bi-variate Taylor monomial.

To ensure the general discrete operator in Eq.~\eqref{eq:Ld} has polynomial consistency of order $p$, the weighting vector \( \mathbf{\Psi}^d \) is evaluated by solving a linear system given by,
\begin{equation}
\mathbf{M}_{i} \cdot \mathbf{\Psi}^{d}{}_{i} = \mathbf{C}^{d}
\end{equation}
where \(\mathbf{M}_i\) is the moment matrix at node \(i\) for recovering the \( p \)-th order consistency, calculated as,
\begin{equation}
\mathbf{M}_{i} = \sum_{j} \mathbf{X}_{ji} \otimes \mathbf{W}_{ji},
\end{equation}
in which \(\mathbf{X}_{ji}\) is the vector of monomials,
\begin{equation}
\label{eq:Xji}
\mathbf{X}_{ji} = \begin{bmatrix} 
x_{ji}, & y_{ji}, & \frac{x_{ji}^2}{2}, & x_{ji}y_{ji}, & \frac{y_{ji}^2}{2}, 
\frac{x_{ji}^3}{6}, & \frac{x_{ji}^2y_{ji}}{2}, & \frac{x_{ji}y_{ji}^2}{2}, & \frac{y_{ji}^3}{6}, 
\frac{x_{ji}^4}{24}, & \cdots & \frac{y_{ji}^p}{p!} 
\end{bmatrix}^T,
\end{equation}

To ensure polynomial consistency of order \(p\), we include the first \(n = (p^2 +3p)/2\) terms in Eq.~\eqref{eq:Xji}. As a result, the leading order of error becomes \(\mathcal{O}( s_{i}^{p+1-l} )\), where $l$ is the order of the derivative being approximated. For example, the LABFM approximation of the first-order derivative converges with \(s^p\), while the second-order derivative with \(s^{p-1} \).

\section{$p$-adaptivity procedures}
\label{sec:p-adapt}
\subsection{Local refinement indicator}
A crucial component of any adaptive strategy is the refinement indicator, the criterion used to decide where the resolution/order should be refined or de-refined. The refinement indicators can be broadly classified into three categories: feature-based, error-based, and goal-oriented\cite{naddei2019comparison}. Feature-based indicators\cite{wackers2014combined}, which identify the key physical features that need to be fully resolved, are computationally cheap and easy to implement, making them highly effective for capturing sharp features like shockwaves; however, their main drawback is that they are not a true measure of numerical error and can fail to refine regions where the error is large but the solution is smooth. In contrast, error-based indicators are more mathematically rigorous, evaluating the error in numerical solutions based on residuals or recovery-based methods\cite{banks2012,ebrahimnejad2015}. They offer the advantage of providing a more reliable and general approximation of the actual discretisation error, but at the cost of increased computational expenses. The primary advantage of goal-oriented (adjoint-based) indicators \cite{fidkowski2011} is their efficiency in optimising the resolution for a specific quantity of interest (e.g., aerodynamic lift); their disadvantage, however, is their high memory requirements and high computational cost for the solution of the adjoint problem\cite{naddei2019comparison}.  

In this study, we utilise the concept of discretisation error-based indicators to provide direct control over the numerical accuracy of the solution. Inspired by the small-scale energy density indicator \cite{kuru2016adaptive,naddei2019comparison} and the spectral decay indicator \cite{persson2006}, which involves projecting the solution onto a reduced-order space for error estimation and has proven effectiveness in $p$-adaptive mesh-based simulations of fluid problems \cite{naddei2019comparison} , this work formulates the indicator in a similar manner. However, in strong-form high-order mesh-free methods, the local construction of lower-order solutions requires substantial additional computational cost, which diminishes the benefits of adaptivity.

To enable a computationally efficient adaptive procedure, the proposed refinement indicator is formulated based on the Laplacian operator, which can be readily integrated into the existing derivative evaluations of a numerical solver. The indicator is also expected to provide a useful measure of the solution error, as the discretisation accuracy of the Laplacian operator plays an important role in determining the numerical accuracy of fluid flow simulations. The dimensionless refinement indicator \(\eta_i\) at node \(i\) is defined as,

\begin{equation}
\eta_i = \left| \widehat{L}^{L}(\phi^*)_{i,h,p} - \widehat{L}^{L}(\phi^*)_{i,h,p -\delta p} \right|
\label{eq:indicator}
\end{equation}
where $\phi^*=\phi/\phi_c$ is the dimensionless field variable, with $\phi_c$ denoting its characteristic scale. The dimensionless Laplacian operator is defined as $\widehat{L}^{L}=L_c^2 L^{L}$, where $L_c$ is the characteristic length scale. The subscripts $h$ and $p$ denote the stencil size and polynomial order, respectively, while $\delta p$ represents the variation of the polynomial order in the refinement. In this work, only even polynomial orders ($p=2,4,6,8,\ldots$) are considered for a monotonic change of error magnitude and \(\delta p = 2\) is applied for all the analysis. The term $\widehat{L}^{L}(\phi^*)_{i,h,p-\delta p}$ represents the dimensionless Laplacian approximation evaluated in the reduced-order approximation space.

\subsection{Refinement criteria}
\label{sec:criteria}
The following refinement criteria are adopted for the refinement or de-refinement of the polynomial order \(p_i\) used for LABFM approximation at a given node \(i\) , after the estimation of the local refinement indicators \(\eta_i\).
\begin{equation}
p_i=\left\{
\begin{array}{ll}
p_i + \delta p & \text{if} \; \eta_i > \varepsilon_U \; \& \; p_i < p_{\text{max}} \\[10pt]
p_i - \delta p & \text{if} \; \eta_i < \varepsilon_L \; \& \; p_i > p_{\text{min}}\\[10pt]
p_i & \text{otherwise} \; 
\end{array}
\right.
\end{equation}
where \(\varepsilon_U\) and \(\varepsilon_L\) are the user-defined limits of the indicator (\(\varepsilon_U >\varepsilon_L\)). Their values depend on the specific requirement of error level and typically range between [\(10^{-8}\), \(10^{-1}\)].  \(p_{max}\) and \(p_{min}\) are respectively the maximum and minimum polynomial order. In this study, \(p_{max}\) is set as $8$, while the \(p_{min}\) is set as $4$.

Since the approximation with a reduced polynomial order employs a smaller stencil compared to that of the original polynomial order, the refinement indicator can be evaluated within the same loop used for the existing Laplacian approximation, incurring only a minor increase in computational cost. The polynomial order distribution generated at the current time step is used as the initial condition for the subsequent time step. As flow features are advected, the refinement indicator may yield large values, prompting an increase in the local polynomial order to accurately capture these features. By assigning lower polynomial orders with smaller stencils in regions of lesser complexity, the overall computational cost is reduced while maintaining the error within a user-defined threshold.

By adding and subtracting the exact dimensionless Laplacian, Eq.~(\ref{eq:indicator}) can be rewritten as
\begin{equation}
\eta_i = \left| [\widehat{L}^{L}(\phi^*)_{i,h,p} - \widehat{\nabla}^2 (\phi^*)_{i}] - [\widehat{L}^{L}(\phi^*)_{i,h,p-\delta p} - \widehat{\nabla}^2 (\phi^*)_{i}]\right|
\end{equation}
where $\widehat{\nabla}^{2}(\phi^*)_i$ denotes the exact dimensionless Laplacian at node $i$. Defining the local Laplacian discretisation error as $\varepsilon^{L}_{i,h,p} = \widehat{L}^{L}(\phi^*)_{i,h,p} -\widehat{\nabla}^{2}(\phi^*)_i$, the indicator becomes
\begin{equation}
\eta_i =
\left|
\varepsilon^{L}_{i,h,p}
-
\varepsilon^{L}_{i,h,p-\delta p}
\right|.
\label{eq:rew_indic}
\end{equation}

This formulation shows that the refinement indicator quantifies the difference between the Laplacian discretisation errors associated with the original and reduced polynomial orders. Since the discretisation error generally decreases with increasing polynomial order for sufficiently smooth solutions, the reduced-order error $\varepsilon^{L}_{i,h,p-\delta p}$ is expected to dominate. Under this assumption, the indicator provides an approximation to the magnitude of the reduced-order Laplacian discretisation error, i.e., $\eta_i \approx |\varepsilon^{L}_{i,h,p-\delta p}|$. The refinement criterion can therefore be used to regulate the estimated error associated with the reduced-order approximation. The effectiveness of this approach is examined in Section~\ref{sec:tests}.

\section{Numerical tests}
\label{sec:tests}
In this section, the $p$-adaptive LABFM is applied to a range of test cases. First, the convergence behaviour of the proposed method is examined to assess how the approximation error is controlled by the thresholds employed in the refinement criteria. Next, the accuracy and performance of the $p$-adaptive method are demonstrated using the travelling wave solution and the periodic solution of the viscous Burgers' equation, as well as the Kelvin–Helmholtz instability. Finally, the proposed method is applied to the simulation of compressible reacting flows in complex geometries, where it is validated against free-propagating hydrogen flames, before being applied to the simulation of flame propagation in porous materials.

\subsection{\label{sec:level2}Convergence test}
This test case evaluates the fundamental error behaviour of gradient and Laplacian approximations using the $p$-adaptive LABFM, over a range of resolutions. A sinusoidal function and a super-Gaussian function are chosen as test functions. The sinusoidal profile is smooth with uniformly distributed errors, making it ideal for analysis of the performance of refinement indicators with varying error thresholds. In contrast, the super-Gaussian features a localised jump with steep gradients, providing a benchmark for testing the adaptive algorithm.

In this case, a unit square domain defined by \((x,y) \in [0,L_c] \times [0,L_c]\), with $L_c=1$, is considered. The tested sinusoidal function is given by
\begin{equation}
\phi(x,y) = \sin\left(k\pi\hat{x}\right)\sin\left(k\pi\hat{y}\right)
\label{eq:sintest}
\end{equation}
where \(\hat{x}=x-0.1453\), \(\hat{y}=y-0.16401\), and \(k\) is the wave number, set as $k=2$ in this case. Note that a pseudo-random offset is introduced into the function to add asymmetry. This prevents the cancellation of numerical errors that might otherwise occur when analysing symmetric function profiles.

For polynomial orders $p=4$, $6$, and $8$, the stencil size ratio \(h/s\) is taken as $1.4$, $1.8$, and $2.3$, respectively. Nodes are initially distributed on a uniform Cartesian grid with resolution \(s\) and then subjected to random perturbations. Following~\cite{king2020high}, the maximum perturbation distance is set as \(\epsilon\) and the resulting perturbation magnitude \(\epsilon/s\) is taken as a measure of distribution irregularity, in which \(\epsilon/s=0.5\) is applied in this study. To ensure full stencil support near boundaries, ghost nodes are generated within a boundary strip of width \(2h\), following the same resolution and distribution irregularity as interior nodes. The field values at the boundary nodes are prescribed using the analytical solution. The $p$-adaptive LABFM method is validated by computing approximations of the gradient and Laplacian of \(\phi\) across multiple resolutions. The initial polynomial order is set to \(p=6\), while local adjustments of the polynomial order are governed by the adaptive algorithm according to the defined error threshold. Errors are quantified using the normalised \(L_2\)-norm. Due to the nearly uniform error distribution in approximating the sinusoidal function, we expect the local error control would be reflected in the global error behaviour.

\begin{figure}[h]
    \centering
    \includegraphics[width=0.49\textwidth]{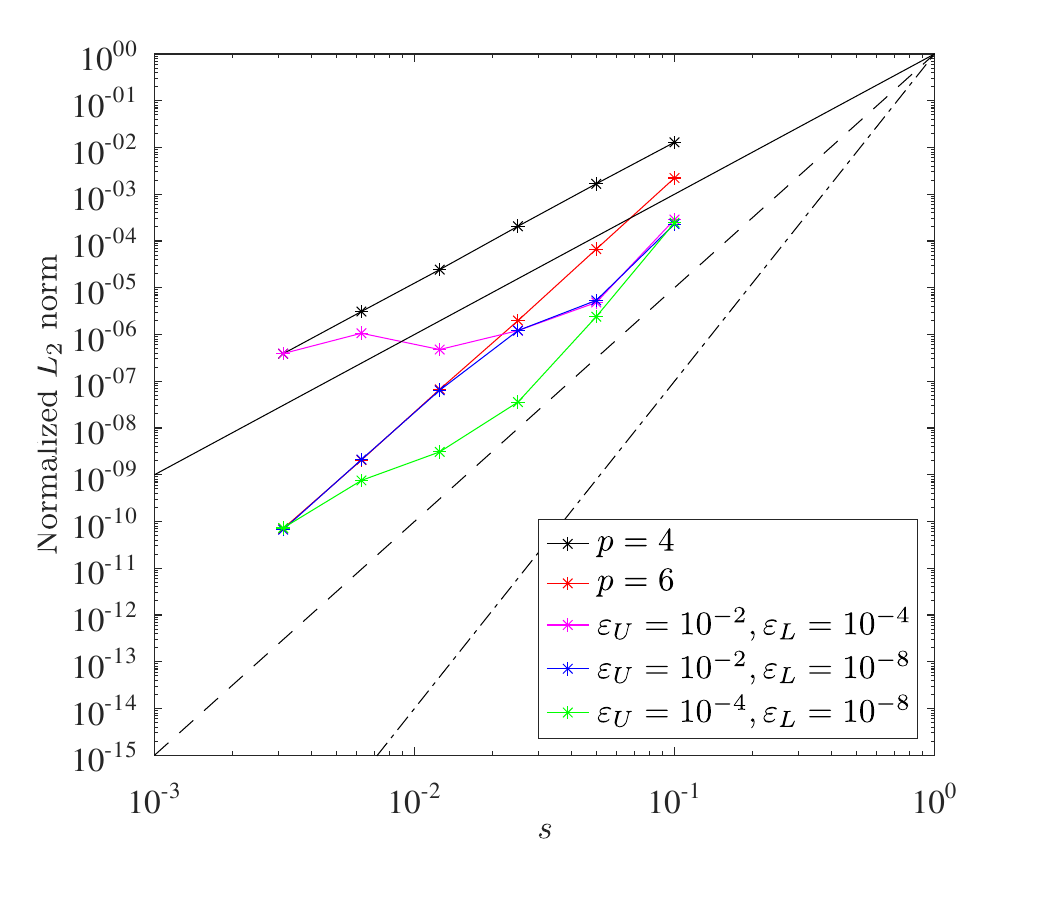}
    \includegraphics[width=0.49\textwidth]{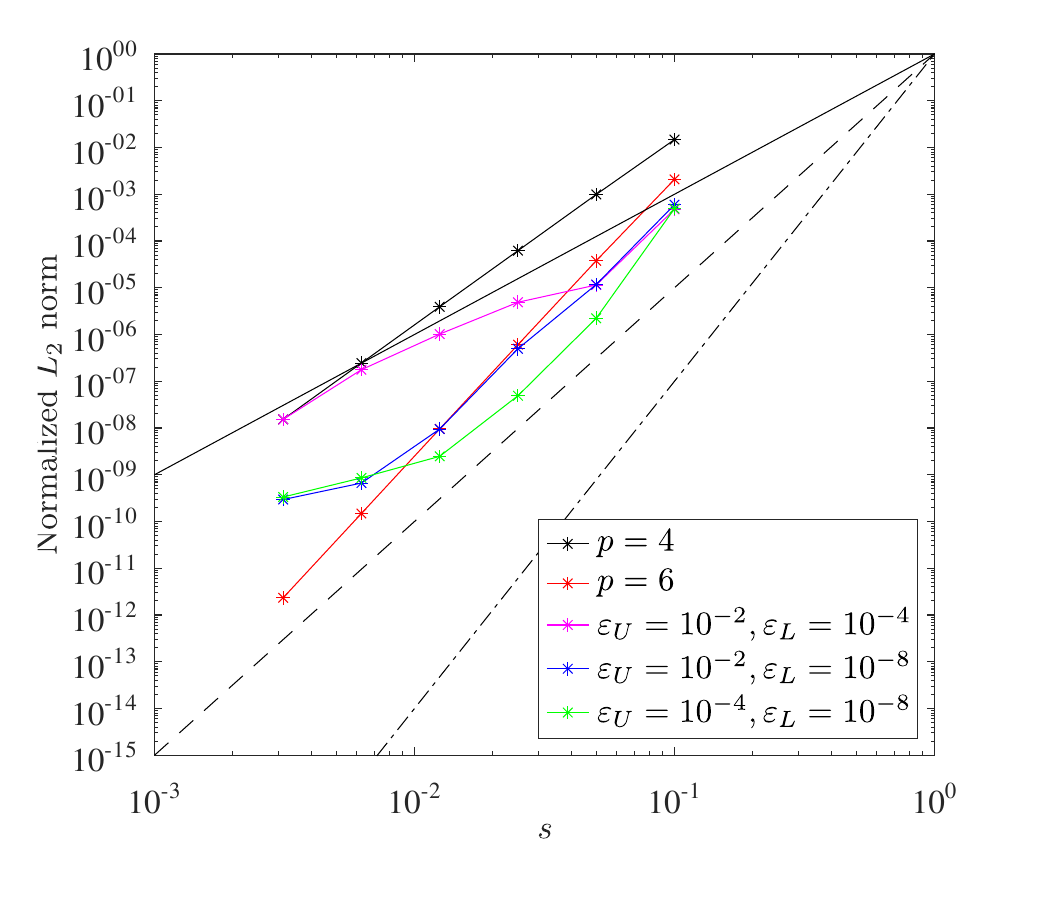}
    \caption{Convergence behaviour of the differential operator approximation using the adaptive and non-adaptive methods, showing Laplacian approximation (left) and gradient approximation (right).\label{fig: convergence_threshold}} 
\end{figure}

Three sets of thresholds are tested, including \(\varepsilon_U=10^{-2}\), \(\varepsilon_L=10^{-4}\); \(\varepsilon_U=10^{-2}\), \(\varepsilon_L=10^{-8}\); and \(\varepsilon_U=10^{-4}\), \(\varepsilon_L=10^{-8}\). The convergence behaviour of the Laplacian and gradient approximation of the sinusoidal function, as $s$ halved successively from $1/10$ to $1/320$, is shown in Fig.~\ref{fig: convergence_threshold}. As discussed in section~\ref{sec:criteria}, the error indicator reflects the error magnitude in the reduced order space. Therefore, in this case, as the resolution varies, a refinement or de-refinement of the polynomial order occurs when the error at that resolution projected to the lower-order space reaches a threshold. For example, this is indicated in the Laplacian approximation by the shift from the convergence rate of $p=6$ to $p=4$ as resolution increases in the line representing \(\varepsilon_U=10^{-2}\), \(\varepsilon_L=10^{-4}\) and the shift from the convergence rate of $p=6$ to $p=8$ as resolution decreases in the line representing \(\varepsilon_U=10^{-4}\), \(\varepsilon_L=10^{-8}\), both of which occur at a resolution of $s=1/40$. In the former case, the error reaches lower threshold \(\varepsilon_L=10^{-4}\) while in the latter it hits the upper threshold \(\varepsilon_U=10^{-4}\). When the error exceeds the upper threshold \(\varepsilon_U\), the interpolation switches to a higher $p$ order for a greater convergence rate. Conversely, if the error is below the lower threshold \(\varepsilon_L\), the interpolation with a lower $p$ order is activated. In between the upper and lower thresholds, the original $p$ order remains (or, globally, a mix of different $p$ orders for non-uniformly distributed errors). This implies that the algorithm will allocate more computation resources to regions identified as having high error, defined by the upper threshold, while less computational focus will be placed on regions with low error, defined by the low threshold, enabling efficient computation according to the predefined accuracy level. It is found that the refinement procedure can also effectively adjust the error behaviour of gradient approximation according to the given thresholds. While the change in convergence rate is observed to occur earlier for the gradient approximation than for the Laplacian approximation as the spatial resolution increases, this is reasonable because the error from the gradient approximation is often smaller than that of the Laplacian, with a convergence rate typically one order of magnitude higher.

In the following test, the super Gaussian function is adopted, given by,
\begin{equation}
\phi(x,y) = \exp\left[-32\pi^4(\hat{x}^4 + \hat{y}^4)\right]
\end{equation}
where the same pseudo-random offset is introduced for $\hat{x}$ and $\hat{y}$ as in the sinusoidal test function, i.e., equation~(\ref{eq:sintest}).

\begin{figure}[h]
    \centering
    \includegraphics[width=0.49\textwidth]{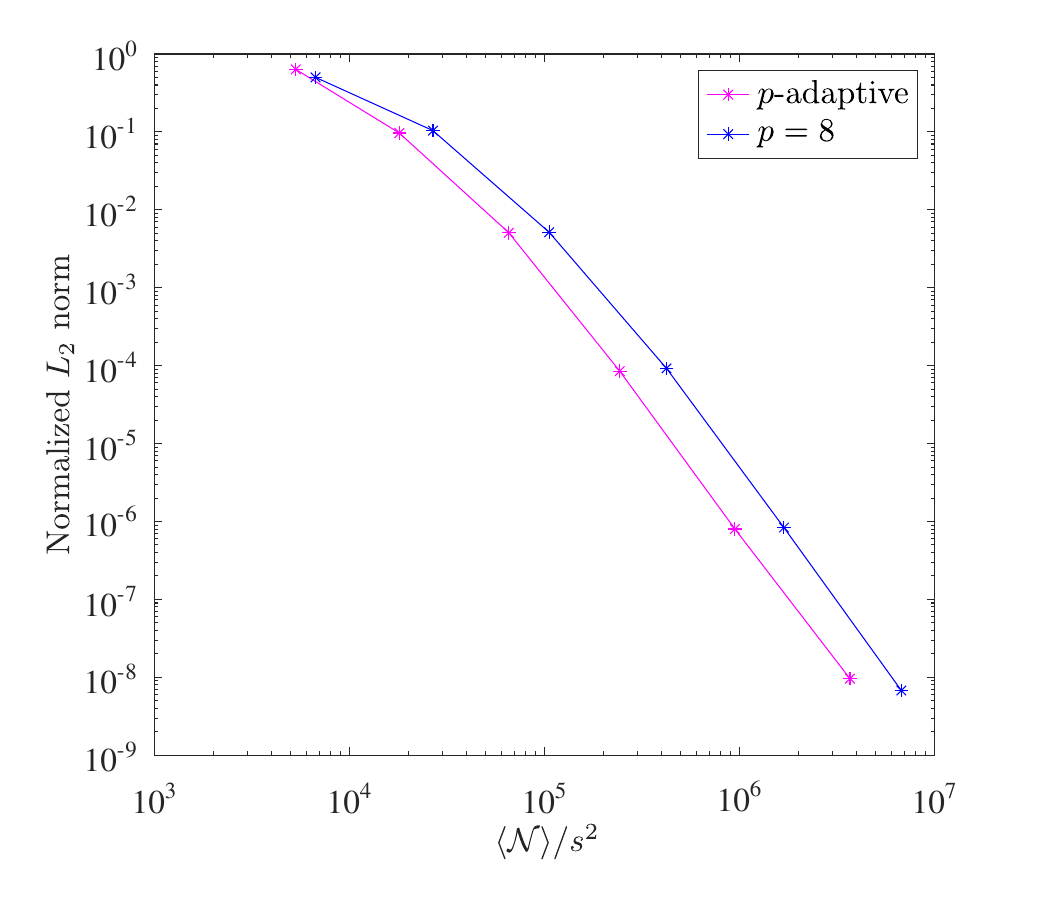}
    \includegraphics[width=0.49\textwidth]{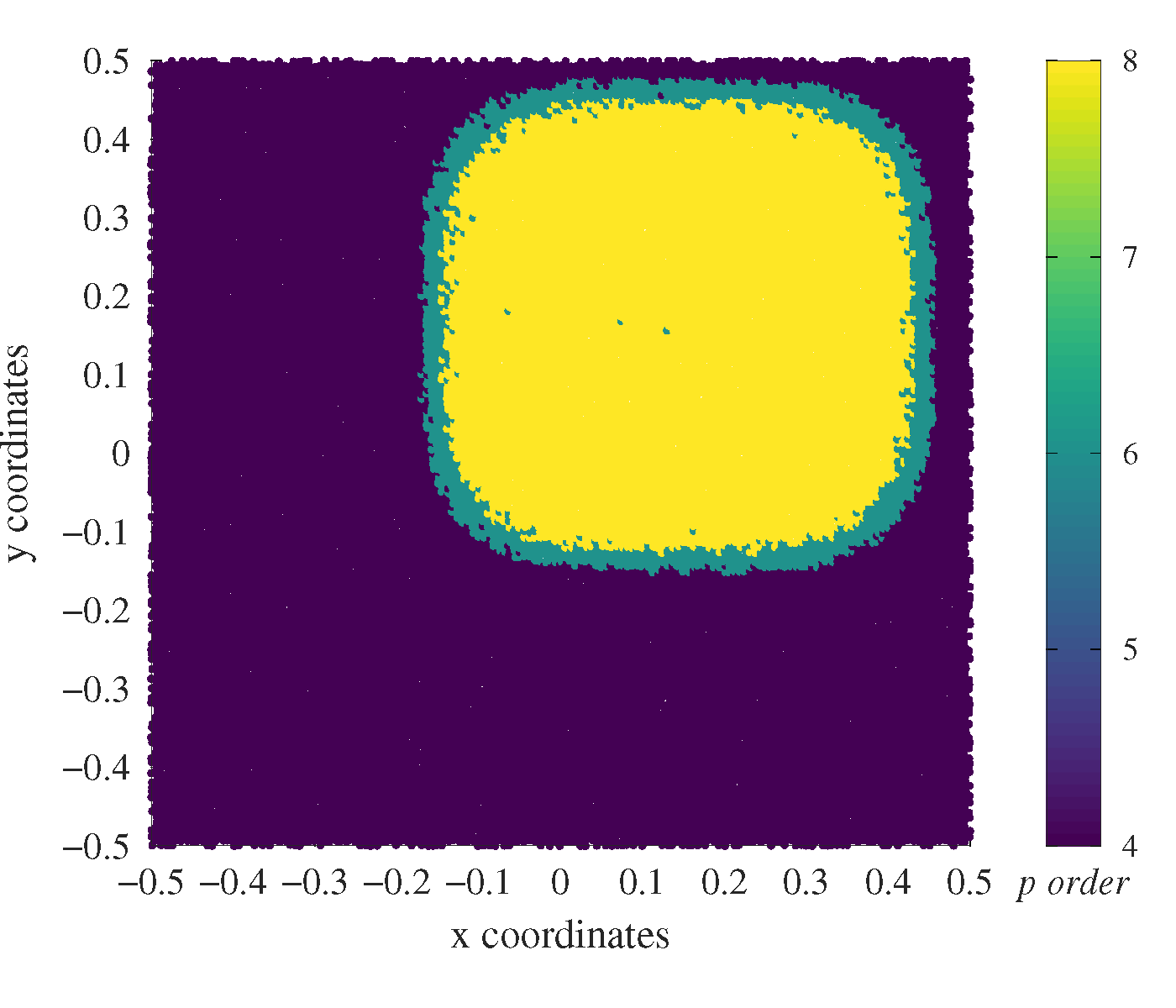}
    \caption{Convergence test results using the adaptive method, showing error versus cost with comparison against the non-adaptive method (left) and local polynomial order $p$ (right). \label{fig: convergence}} 
\end{figure}

Fig.~\ref{fig: convergence} shows the relationships between the error and the cost by $p=8$ and $p$-adaptivity [\(\varepsilon_U=10^{-5}\), \(\varepsilon_L=10^{-8}\)], as well as the local polynomial order when using the adaptive method. Herein, the cost is quantified by a measure, the average neighbour count divided by the resolution $\left\langle\mathcal{N}\right\rangle/s^{2}$, which is proportional to the computational cost. As shown in Fig.~\ref{fig: convergence}, it is clear that the adaptive algorithm effectively captures regions with steep gradients that require high-order polynomial reconstruction.  For the same computational cost, the $p$-adaptive algorithm achieves higher accuracy than a fixed polynomial order. For example, for $\left\langle\mathcal{N}\right\rangle/s^{2}=10^6$, the Laplacian approximation with \(p = 8\) leads to the normalised $L_2$ norm error of the order of \(10^{-5}\), while the $p$-adaptivity results in an error of the order of \(10^{-6}\).

\subsection{Viscous Burgers’ equation}
In this case, the $p$-adaptive LABFM is tested by solving the viscous Burgers’ equation, which is a nonlinear advection--diffusion equation that develops steep viscous layers.  The viscous Burgers' equation is written as,
\begin{equation}
    \frac{\partial \mathbf{u}}{\partial t} + \mathbf{u} \cdot \nabla \mathbf{u} = \frac{1}{Re} \nabla^2 \mathbf{u}
\label{eq:vb}
\end{equation}
where \(Re\) is the Reynolds number. We use LABFM (adaptive or non-adaptive) to evaluate spatial derivatives, and an explicit fourth order Runge-Kutta (RK4) scheme for time integration. The adopted time step (dimensionless)  for the solution of Eq.~(\ref{eq:vb}) is given by,
\begin{equation}
    \delta t = \min\left(\frac{0.2h}{\max\left\lVert\mathbf{u} \right\lVert}, 0.05h^2Re\right)
\end{equation}
\subsubsection{Travelling wave test case}
The first test involves a travelling wave moving diagonally along the xy-plane, where the steepening caused by non-linear advection is balanced by viscous diffusion. As the Reynolds number increases, the wave becomes steeper. The computational domain is a unit-square defined by \((x,y) \in [0,L_c] \times [0,L_c]\), with $L_c=1$, and discretised using a perturbed Cartesian node distribution with \(\epsilon/s=0.2\). To ensure complete stencil support at boundaries, the node distributions near boundary are extended outward, and the analytical solutions are imposed on these exterior nodes. The initial conditions and analytic solution are given by:
\begin{align}
u(x,y,t) &= \frac{3}{4} - \frac{1}{4\left(1 + e^{Re(-t-4x + 4y)/32}\right)},  \\
v(x,y,t) &= \frac{3}{4} + \frac{1}{4\left(1 + e^{Re(-t-4x + 4y)/32}\right)}. 
\end{align}

We vary the (spatially uniform) resolution from \(s=1/10\) to \(s=1/160\). The initial polynomial order is taken as \(p=8\). The Reynolds number \(Re\) is taken as $200$ and $500$ for two distinct wave profiles. The upper and lower thresholds are set as \(\varepsilon_U=10^{-3}\) and \(\varepsilon_L=10^{-6}\).

\begin{figure}[h]
    \centering
    \includegraphics[width=0.49\textwidth]{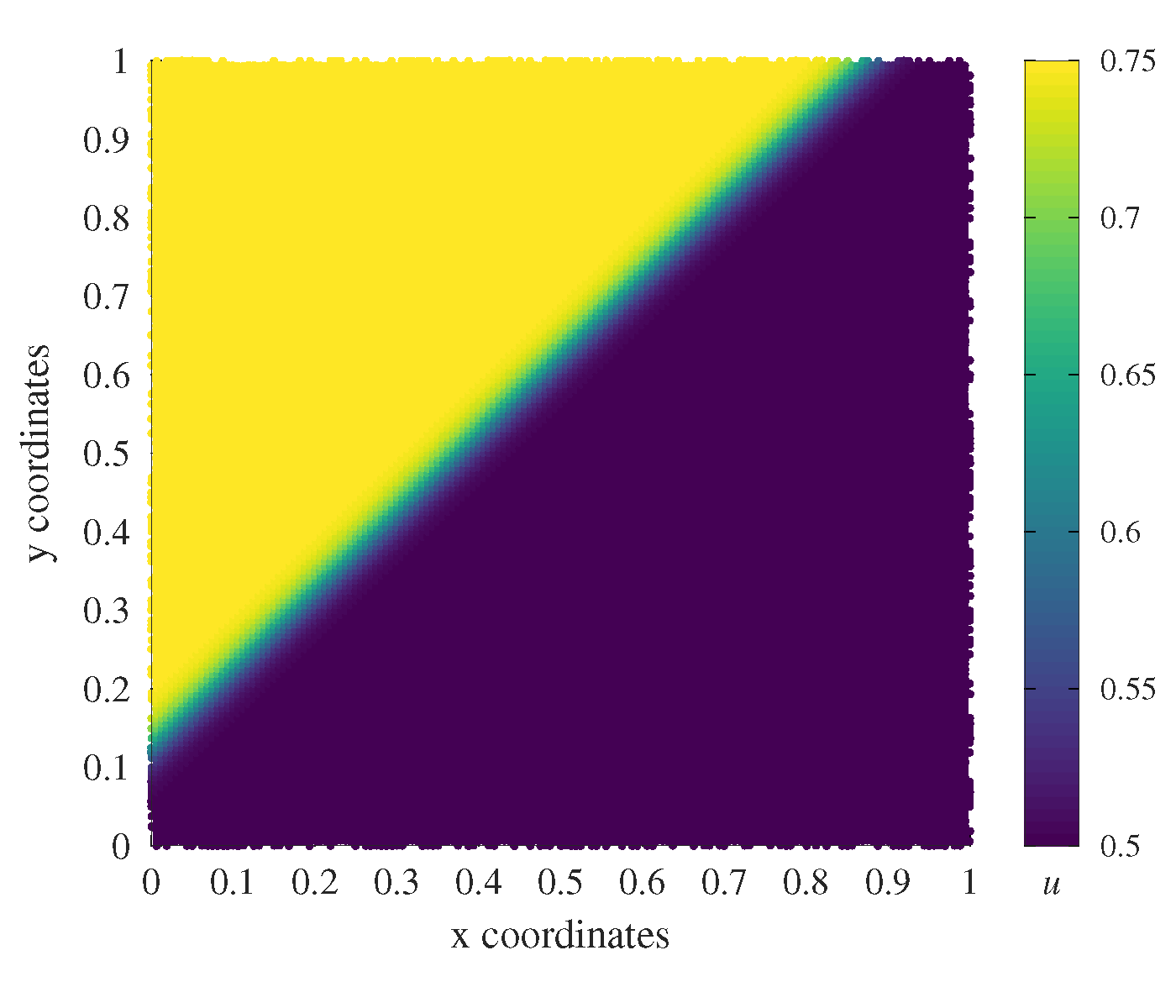}
    \includegraphics[width=0.49\textwidth]{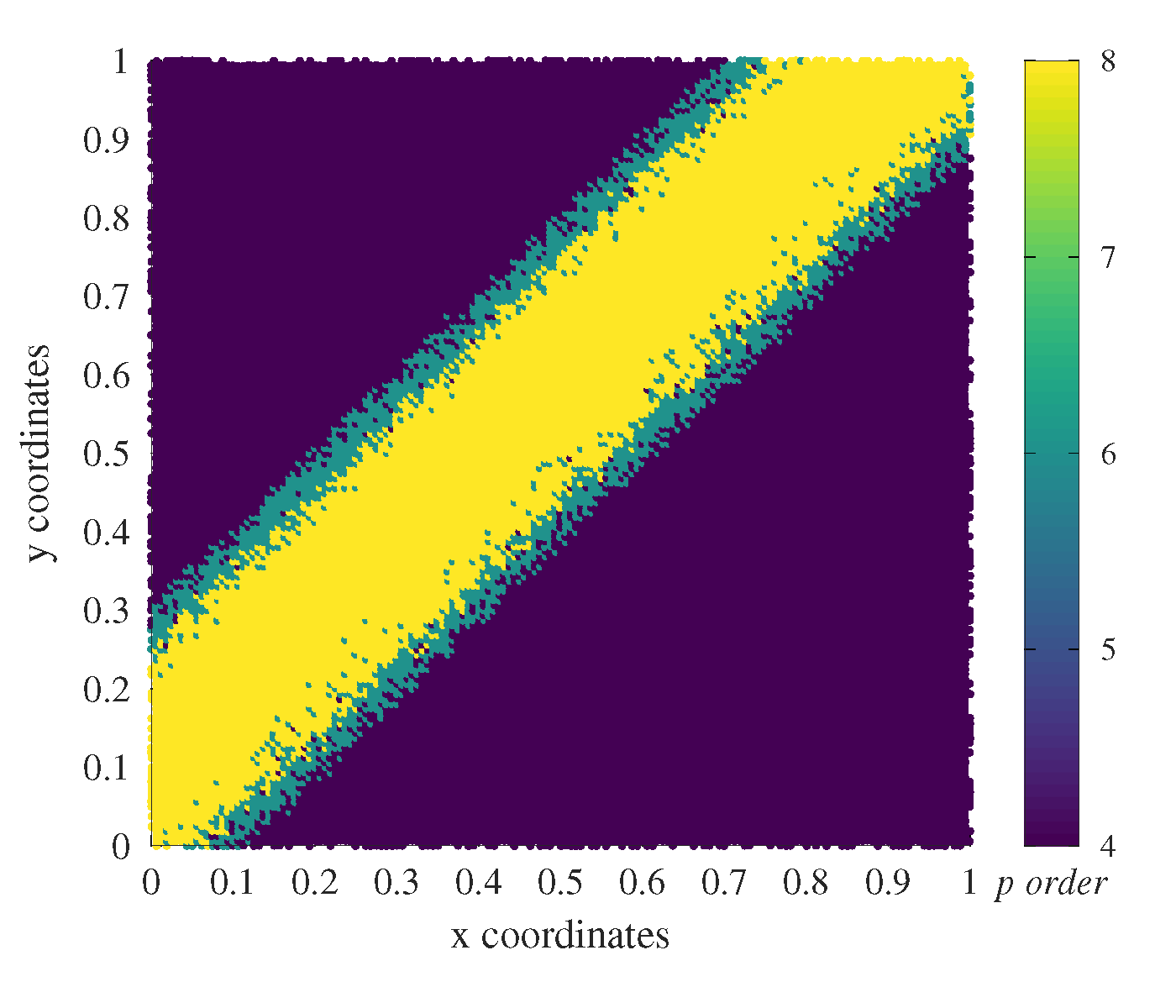}
    \caption{Numerical solution of the travelling wave case using the adaptive method [\(Re = 500\), \(s=1/160\), \(t=0.5\)], showing the value of $u$ (left) and local polynomial order $p$ (right).\label{fig: travelwave_Re=500}} 
\end{figure}

\begin{figure}[h]
    \centering
    \includegraphics[width=0.49\textwidth]{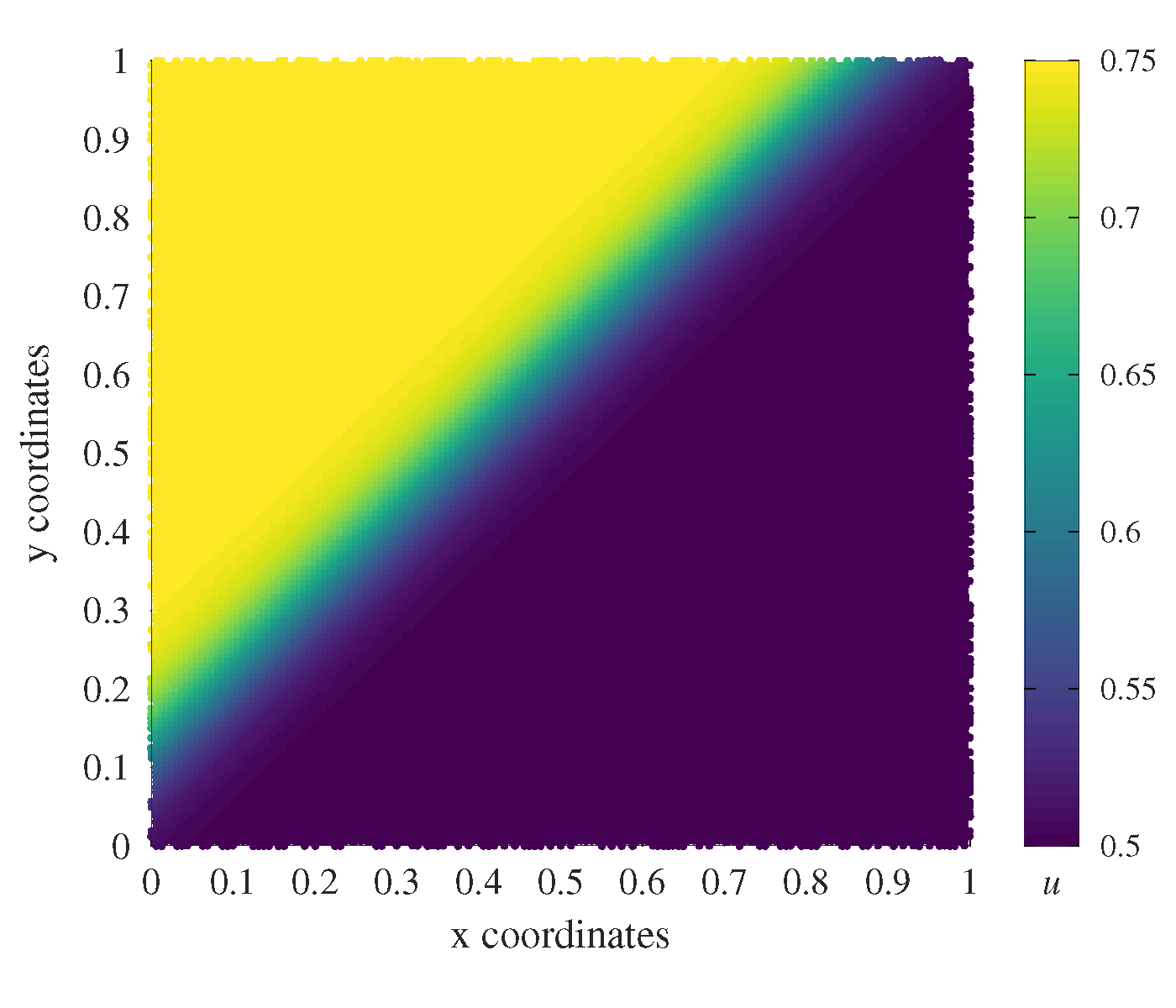}
    \includegraphics[width=0.49\textwidth]{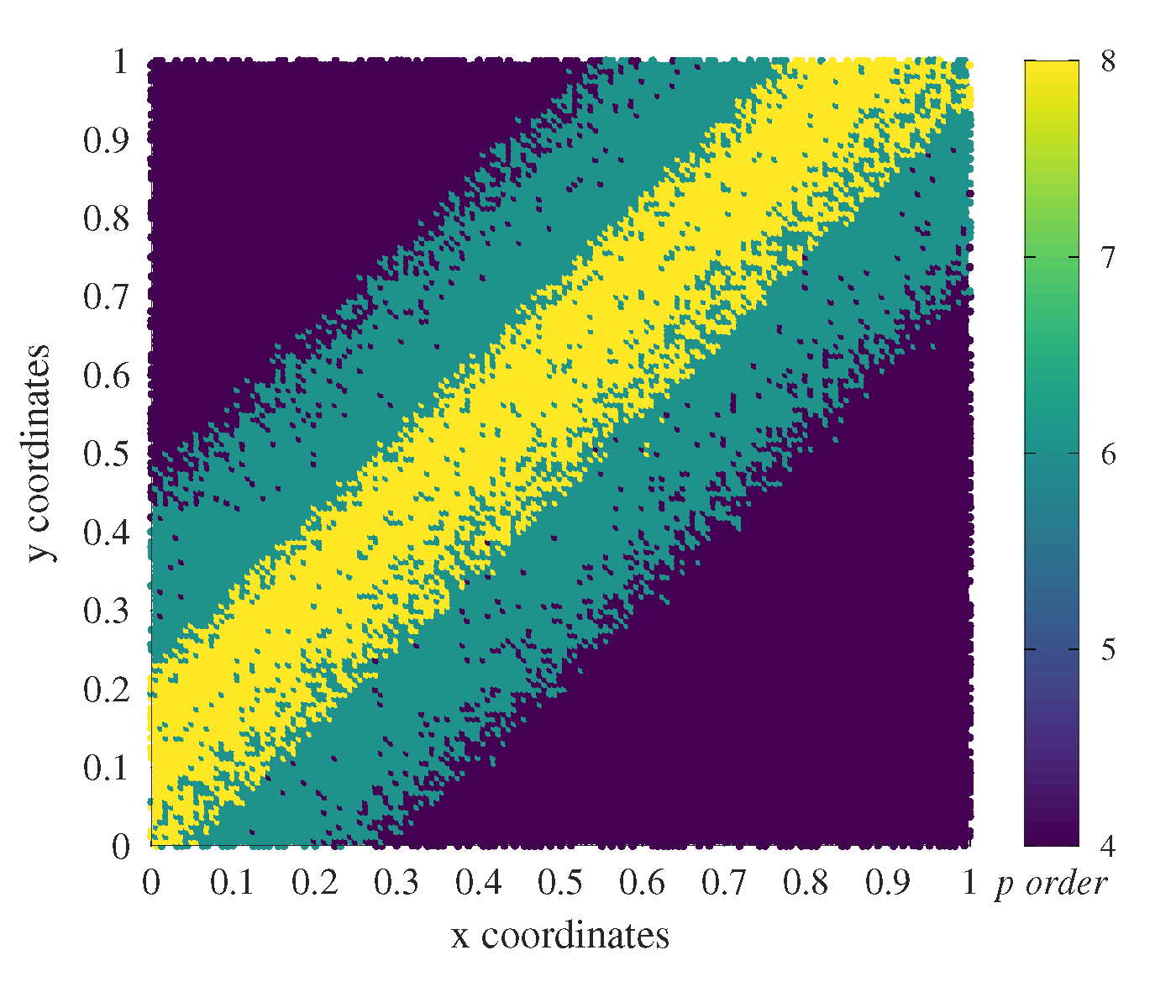}
    \caption{Numerical solution of the travelling wave case using the adaptive method [\(Re = 200\), \(s=1/160\), \(t=0.5\)], showing the value of $u$ (left) and local polynomial order $p$ (right).\label{fig: travelwave_Re=200}} 
\end{figure}

\begin{figure}[h]
    \centering
    \includegraphics[width=0.49\textwidth]{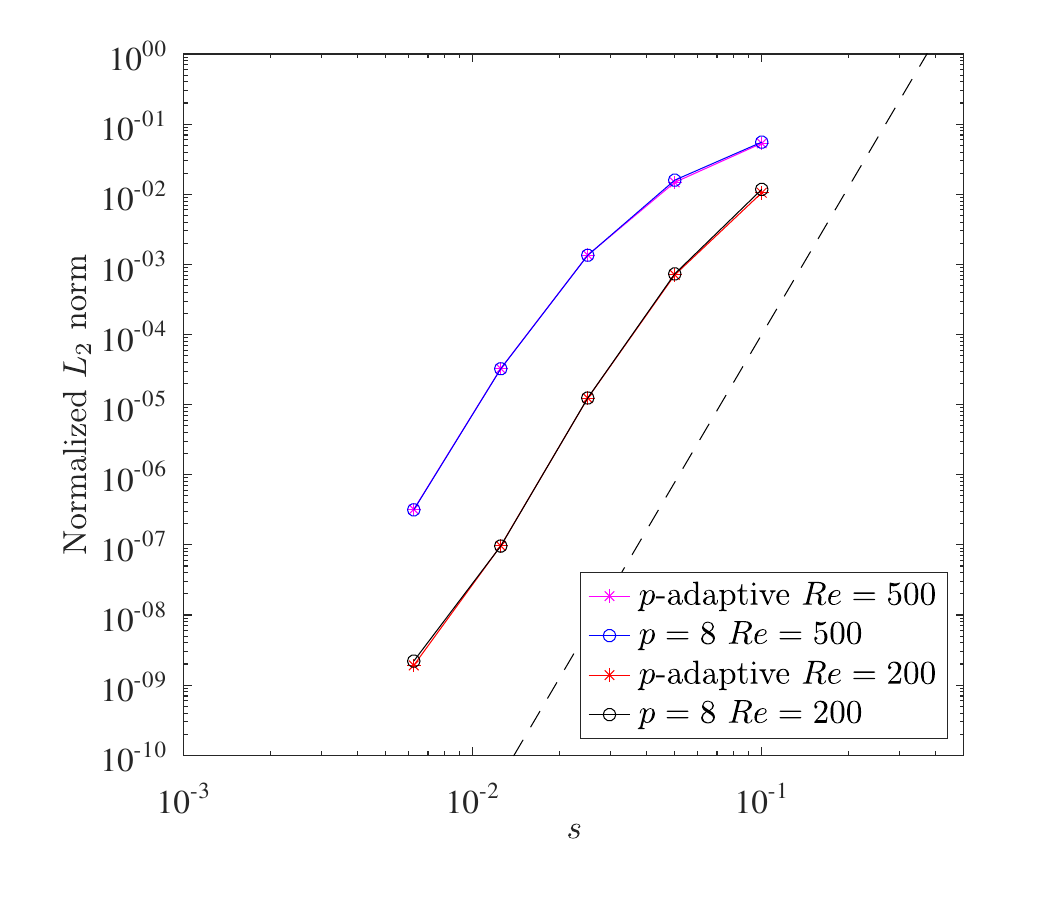}
    \includegraphics[width=0.49\textwidth]{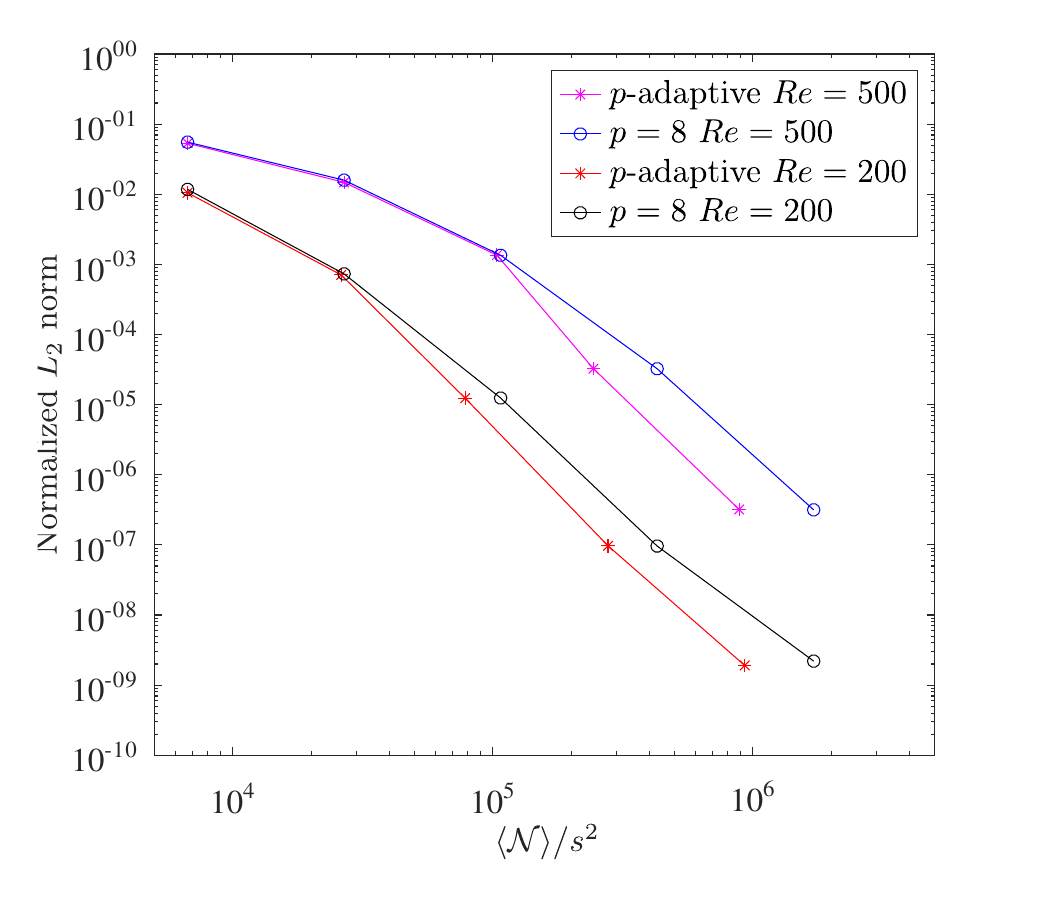}
    \caption{Error analysis of the travelling wave test at $t=1$ for $Re=200$ and $Re=500$, using the $p$-adaptive method and the non-adaptive method with $p=8$. The left panel shows the convergence behaviour of the error, while the right panel shows the error versus computational cost.\label{fig: convergence_travel_wave}} 
\end{figure}

Fig.~\ref{fig: travelwave_Re=200} and Fig.~\ref{fig: travelwave_Re=500} show the numerical solution of $u$ and the local polynomial order for the simulation with $Re=200$ and $Re=500$, respectively, with a node spacing of \(s=1/160\) at \(t = 0.5\). The adaptive algorithm is shown to capture the moving wave profile and effectively locate the region with nonlinearity. It is also noted that the case with $Re=200$ exhibits a narrower range of the highest polynomial order and a wider range of medium polynomial order compared to the case with $Re=500$. This is expected because the wave profile for $Re=200$  is smoother than that for $Re=500$, and the spatial distribution of error remains flatter for $Re=200$. The convergence of error for this case is shown in Fig.~\ref{fig: convergence_travel_wave}. The convergence behaviour of the error by $p$-adaptivity for both $Re=200$ and $Re=500$, across the given resolutions from \(s=1/10\) to \(s=1/160\), is similar to that with a fixed polynomial order of $p=8$.  However, the $p$-adaptive method exhibits a clear reduction in computational cost, as shown in the right panel of Fig.~\ref{fig: convergence_travel_wave}. This reduction is more pronounced at higher Reynolds numbers, where the wave profile becomes steeper and narrower. For the node spacing of \(s=1/160\), the adaptive algorithm can improve the accuracy by one order of magnitude compared to the non-adaptive method at the same $\left\langle\mathcal{N}\right\rangle/s^{2}$. No cost reduction is observed at the first few coarsest node resolutions. This is attributed to the local error being greater than the upper threshold \(\varepsilon_U=10^{-3}\), resulting in the highest polynomial order activated across the entire domain without any local adaptation. 

\subsubsection{Periodic test case}
The second case refers to the solution of the viscous Burgers' equation on a unit periodic square domain. The initial velocity condition is defined by a sinusoidal function given by \(u(x,y,0) = \sin\left({2\pi x}\right)\) and \(v(x,y,t) = 0\). The analytic solution, derived via the Cole-Hopf transformation was first given in~\cite{cole_1951}. Three initial node spacings \(s=1/40\), \(s=1/80\), and \(s=1/160\) are tested on an unstructured node distribution with the relative perturbation amplitude of \(\epsilon/s=0.2\). The Reynolds number \(Re\) is taken as $100$. The refinement thresholds match those of the travelling wave case, with \(\varepsilon_U=10^{-3}\) for refinement and \(\varepsilon_L=10^{-6}\) for coarsening.

\begin{figure}[h]
    \centering
    \includegraphics[width=0.49\textwidth]{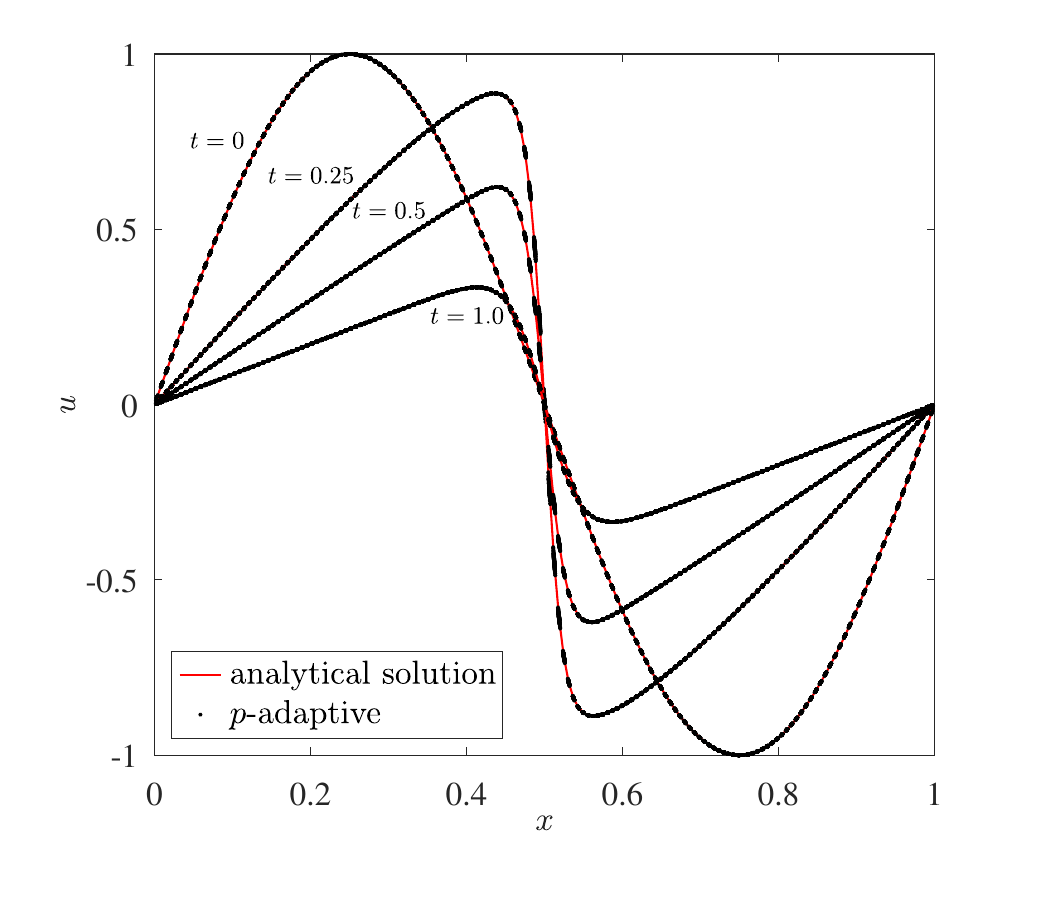}
    \includegraphics[width=0.49\textwidth]{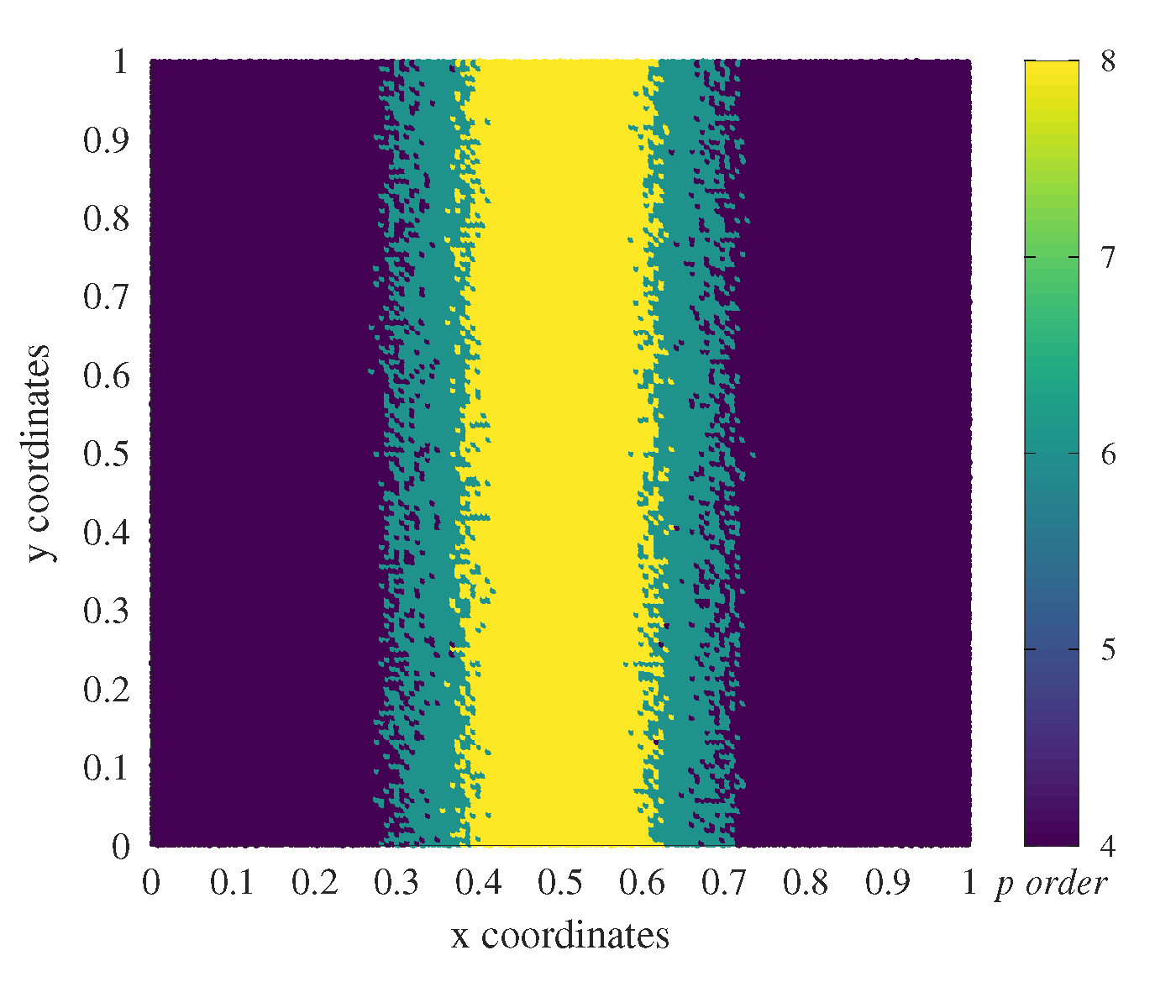}
    \caption{Numerical solution of the periodic viscous Burgers' equation, showing the time evolution of $u$ against analytical solution (left) and local polynomial order $p$ (right) [\(Re = 100\), \(s=1/160\), \(t=0.5\) ].\label{fig: periodic_burgers_solutions}} 
\end{figure}

In Fig.~\ref{fig: periodic_burgers_solutions}, the comparisons of numerical solution against the analytical solution at different time instants are presented, together with the plot of local polynomial order at $t=0.5$. The initial sinusoidal profile gradually develops into a shock due to the contribution of the advection term and then dissipates under the effect of the viscous term. The $p$-adaptive LABFM matches well with the analytical solution. The proposed $p$-adaptive algorithm effectively captures the shock region while allowing lower polynomial orders elsewhere, resulting in a stable, accurate, and efficient solution.

\begin{figure}[h]
    \centering
    \includegraphics[width=0.49\textwidth]{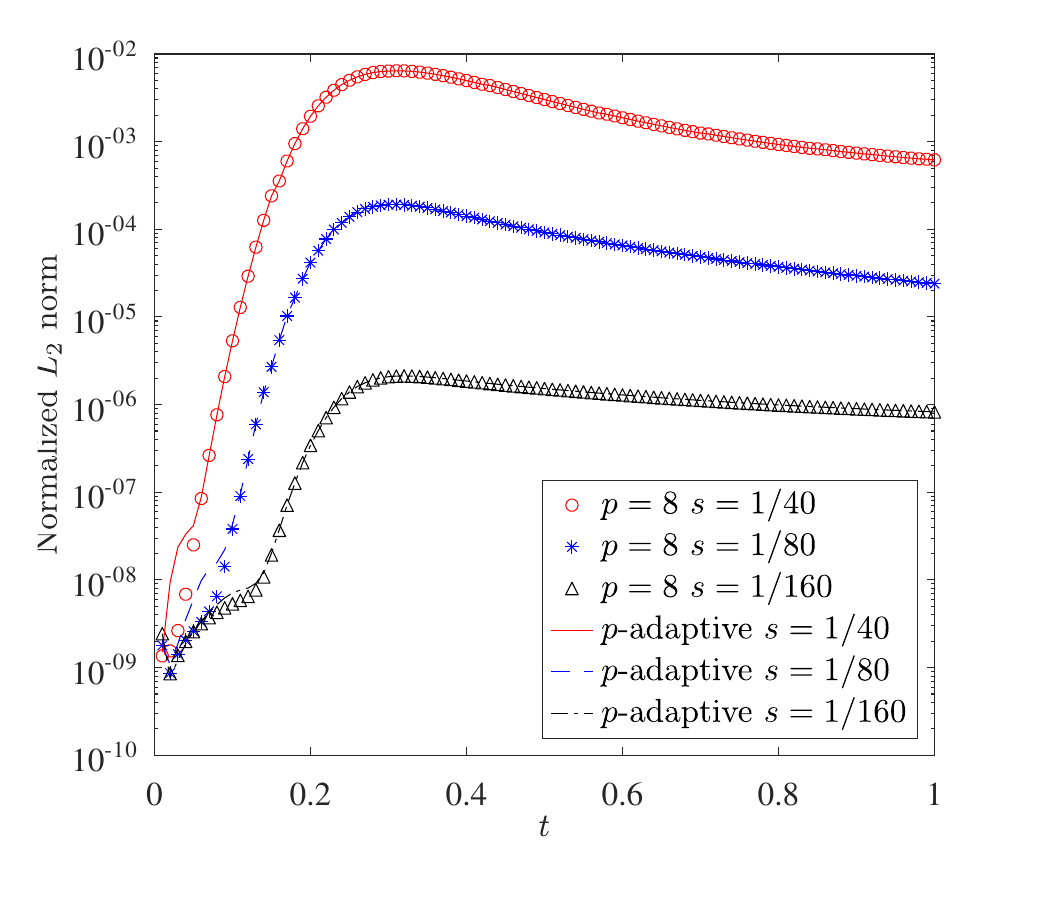}
    \includegraphics[width=0.49\textwidth]{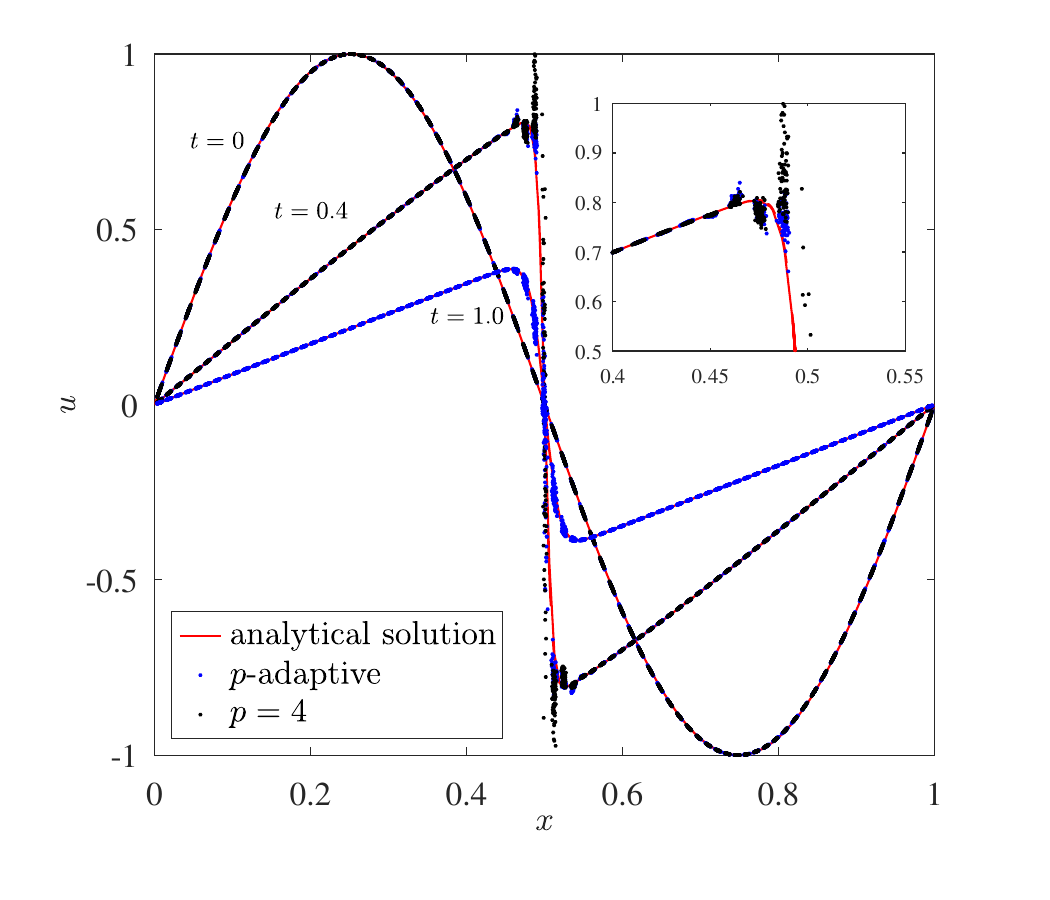}
    \caption{Time evolution of the error in the periodic viscous Burgers' test by the $p$-adaptive method and the non-adaptive method with $p=8$ at resolutions of \(s=1/40\), \(s=1/80\), and \(s=1/160\), using \(\epsilon/s=0.2\) and \(Re = 100\) (left) and time evolution of $u$ against the analytical solution, with the numerical solution obtained by the $p$-adaptive method (blue dots) and the non-adaptive method with $p=4$ (black dots) at resolutions of  \(s=1/80\), using \(\epsilon/s=0.4\) and \(Re = 250\) (right).\label{fig: periodic_burgers_convergence}} 
\end{figure}

The variation of \(L_2\)-norm error with time for the resolutions of \(s=1/40\), \(s=1/80\), and \(s=1/160\) is shown in Fig.~\ref{fig: periodic_burgers_convergence}. Errors increase during the initial advection-dominated phase as the sinusoidal profile steepens and subsequently decrease during the viscosity-dominated phase. Except for the initial few time instants, where the error is very low and reaches the threshold for coarsening, the error from the $p$-adaptivity is nearly identical to that from a fixed polynomial order of 8.

A further investigation is conducted with a configuration where the Reynolds number is increased to $Re = 250$ and the disorder in the node distribution is increased to $\epsilon/s = 0.4$. As shown in the right panel of Fig.~\ref{fig: periodic_burgers_convergence}, the results using a uniform $p = 4$ exhibit overshoots around the shock, which eventually lead to numerical instability and failure of the simulation (no results after around $t = 0.5$). This occurs because the shock profile is under-resolved when using low polynomial orders, causing the solution to become unstable and necessitating higher-order approximations in these regions~\cite{king2020high}. By contrast, the proposed adaptive algorithm, which varies the $p$  order from $4$ to $8$, effectively identifies the shock region and mitigates the instability, leading to an accurate and stable simulation.

\subsection{Kelvin–Helmholtz instability}

Next we study the Kelvin–Helmholtz instability, a fluid instability that occurs at the interface between two fluids due to shear~\cite{kordyban1970,lin1986,price2008}. In this case, an initial sinusoidal perturbation was given to the two-fluid system at the interface. The node distribution without perturbations is used for comparable simulations, as randomly distributed node positions can alter the instability pattern from case to case.

The dimensionless two-phase Navier-Stokes equations for this case are expressed as
\begin{align*}
\frac{\partial \rho}{\partial t} + \bm{u} \cdot \nabla \rho &= -\rho \nabla \cdot \bm{u} \\
\frac{\partial \bm{u}}{\partial t} + \bm{u} \cdot \nabla \bm{u} &= -\nabla p_f + \frac{1}{Re} \nabla^2 \bm{u} \\
\frac{\partial Y}{\partial t} + \bm{u} \cdot \nabla Y &= \frac{1}{Re \, Sc} \nabla^2 Y
\end{align*}
where $p_f$ is the fluid pressure, $Re$ is the Reynolds number, $Fr$ is the Froude number, $Sc$ denotes the Schmidt number, and $At$ represents the Atwood number.

The governing equations are closed by a barotropic equation of state for two fluids, given by,
\begin{equation}
p_f = \frac{1}{Ma^2} \left( \rho - Y - (1 - Y) \frac{1 - At}{1 + At} \right)
\end{equation}

The initial conditions are set as follows:
\begin{align*}
u = 1/2, \quad v = 0, \quad Y = 0, \quad \rho = (1-At)/(1+At) \quad \forall \quad y \in [0.5 + \delta \sin(2\pi x), 1] \\
u = -1/2, \quad v = 0, \quad Y = 1, \quad \rho = 1 \quad \forall \quad y \in [0,0.5 + \delta \sin(2\pi x)]
\end{align*}
where \(\delta\) is the disturbance coefficient for creating the instability, which is taken as 0.05 in this work. 

The time steps (dimensionless) for numerical solution are constrained by the speed of sound, advection and diffusion, given by,
\begin{equation}
    \delta t = \min\left\{ \frac{0.1h}{u_{\text{max}} \left(1 + \frac{1}{Ma} \right)}, \; 0.05 h^2 Re\min\left(1,Sc\right) \right\}
\end{equation}

The flow properties for the two-phase flow are set as \(Re=500\), \(Sc=8\), \(Ma=0.1\), and \(At=0.2\). The node resolution is set as \(s/H=160\). The upper and lower thresholds are taken as \(\varepsilon_U=10^{-1}\) and \(\varepsilon_L=10^{-4}\).

\begin{figure}[h]
    \centering
    \includegraphics[width=0.49\textwidth]{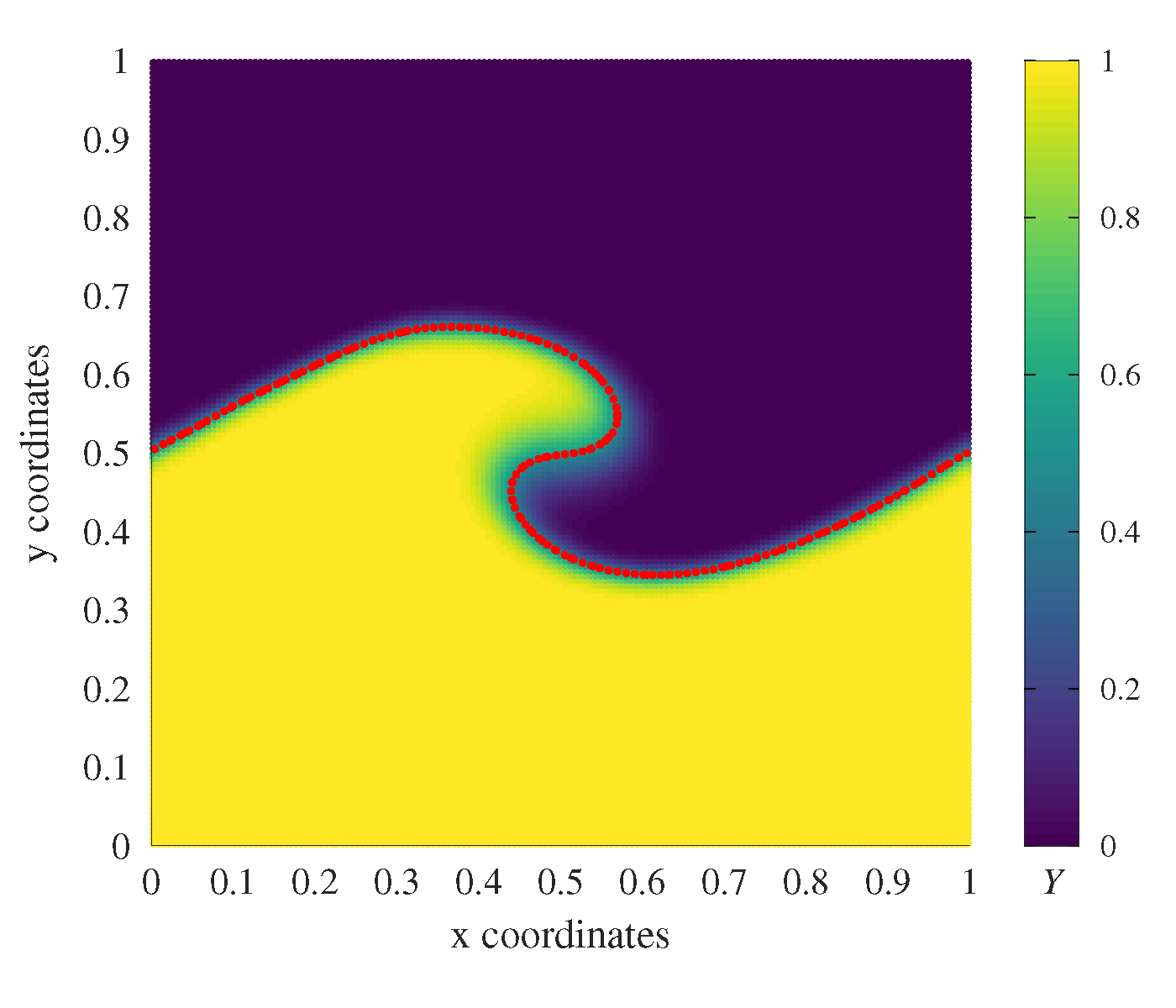}
    \includegraphics[width=0.49\textwidth]{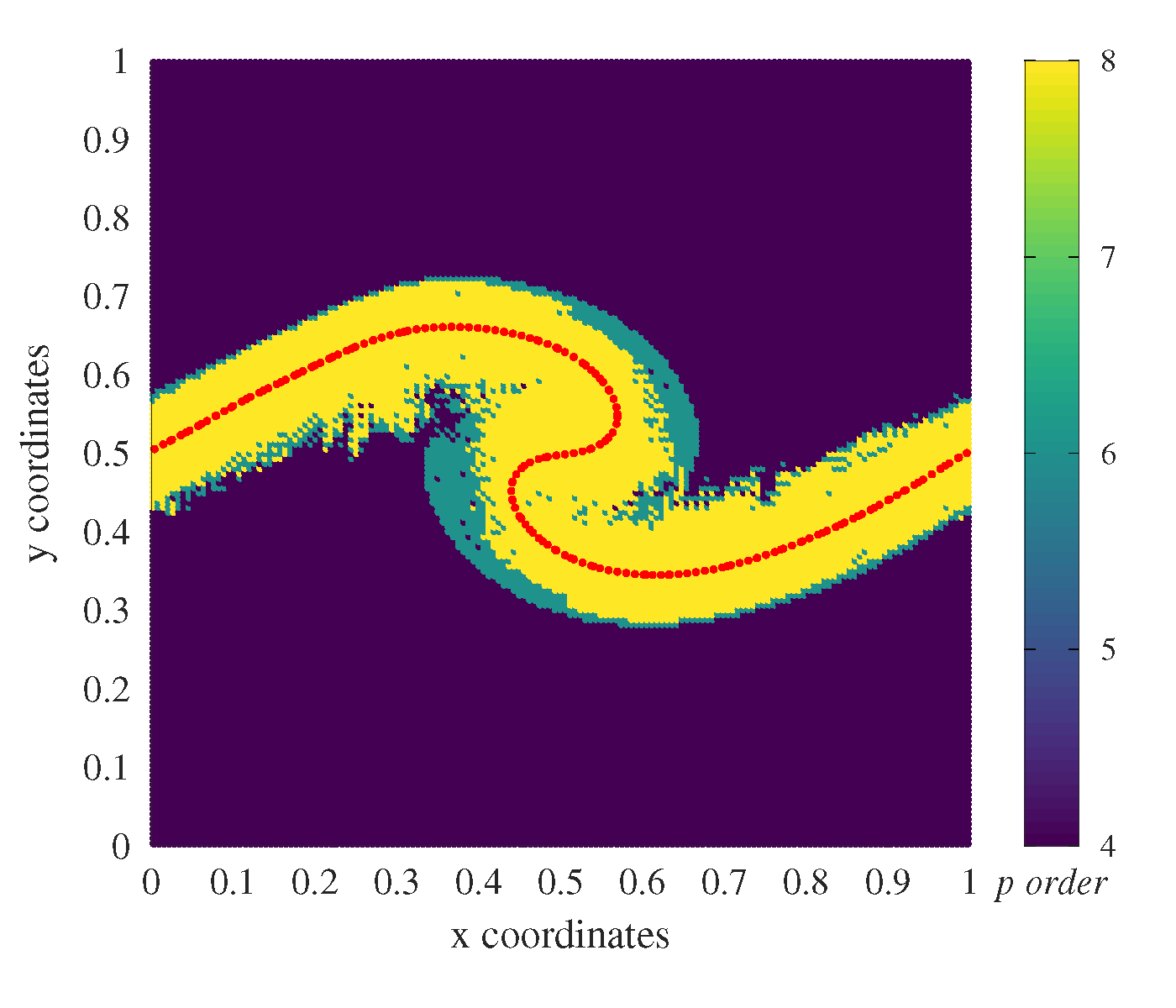}
    \caption{The numerical solution of the Kelvin–Helmholtz instability at $t=1$ using the adaptive method, showing the volume fraction evolution (left) and local polynomial order $p$ (right), where the red line represents the interpolated interface.\label{fig: KH_t=1}} 
\end{figure}

\begin{figure}[h]
    \centering
    \includegraphics[width=0.49\textwidth]{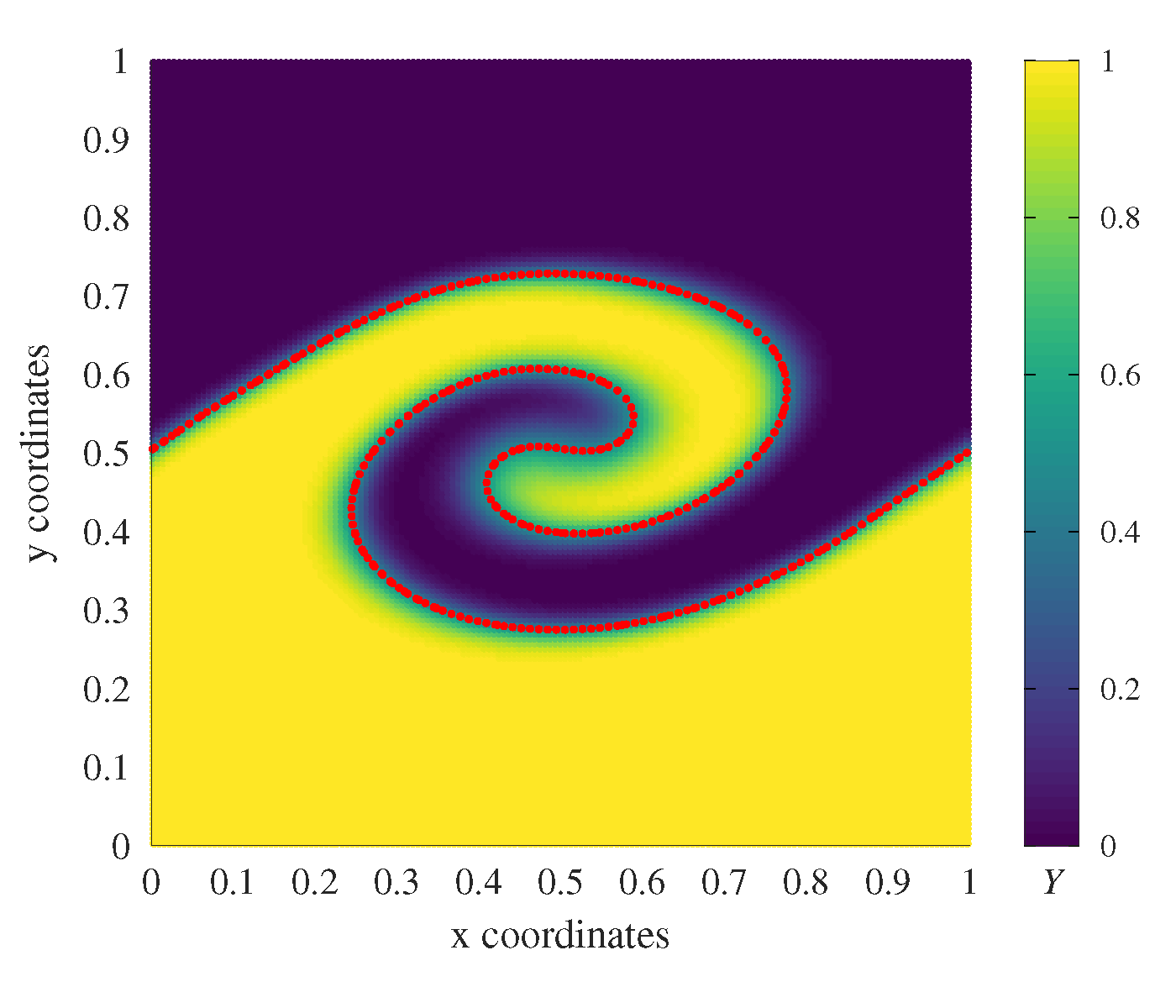}
    \includegraphics[width=0.49\textwidth]{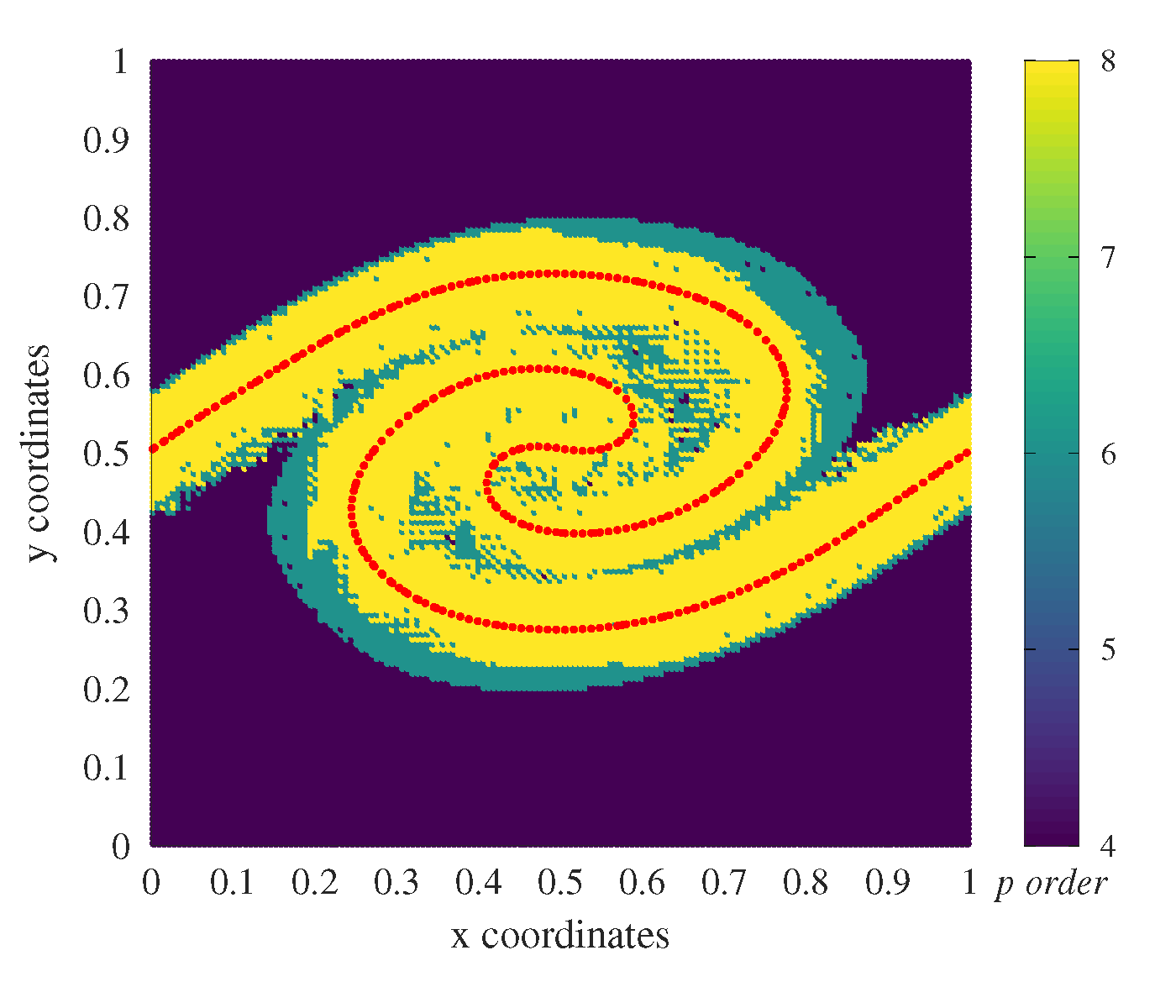}
    \caption{The numerical solution of the Kelvin–Helmholtz instability at $t=2$ using the adaptive method, showing the volume fraction evolution (left), and local polynomial order $p$ (right), where the red line represents the interpolated interface.\label{fig: KH_t=2}} 
\end{figure}

Figures~\ref{fig: KH_t=1} and ~\ref{fig: KH_t=2} show the evolution of the volume fraction at $t = 1$ and $t= 2$ , respectively, along with their corresponding local polynomial orders. The interpolated interface of the two-phase flow (defined as the location where \(Y=0.5\) ) is also plotted. The growth and roll-up of perturbations into the formation of typical cat-eye structures in the Kelvin–Helmholtz instability~\cite{chen2011,shadloo2011} is reproduced using the $p$-adaptive LABFM and the evolution of local polynomial orders follows the changing vortex structures.

\begin{figure}
    \centering
    \includegraphics[width=0.49\textwidth]{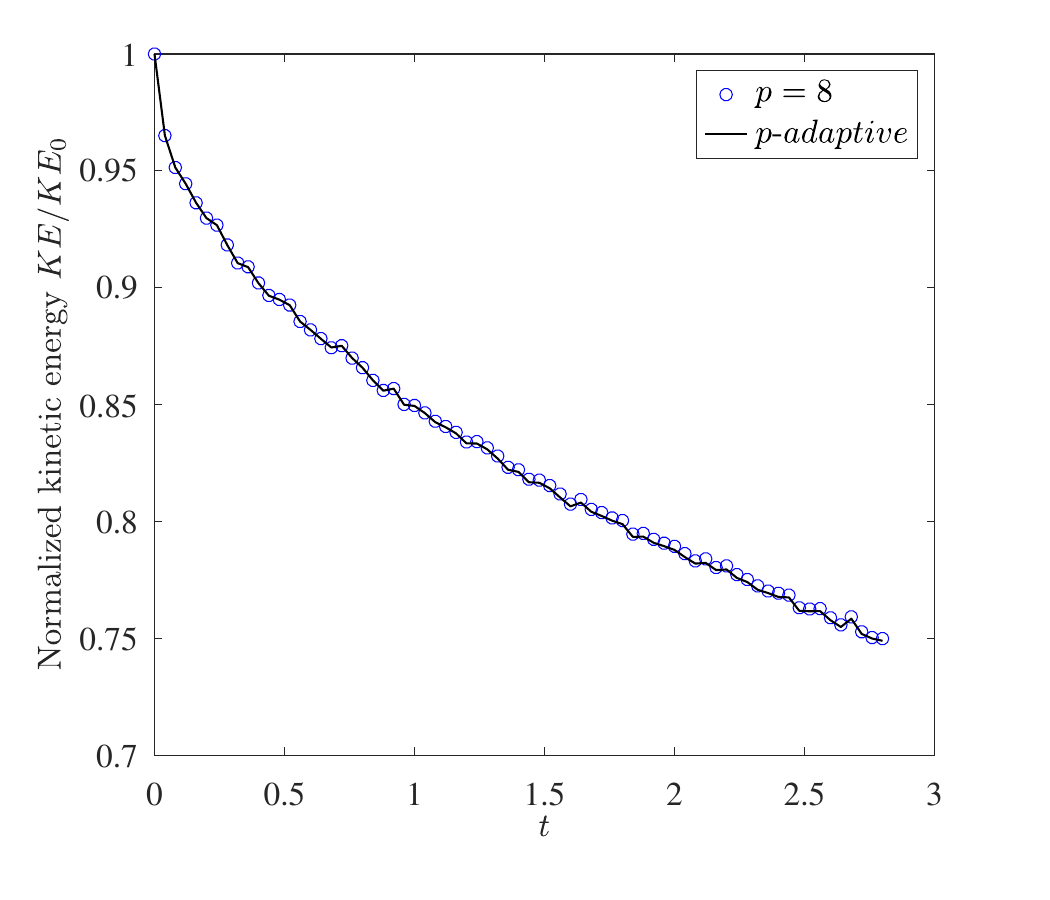}
    \caption{Time evolution of the normalised kinetic energy $KE/KE_0$ for the non-adaptive method with $p=8$ and the adaptive method.}\label{fig: KH_KE}
\end{figure}

Fig.~\ref{fig: KH_KE} presents the calculated time evolution of the total kinetic energy ($KE$), normalised by the initial kinetic energy ($KE_0$), using a fixed polynomial order ($p = 8$) and $p$-adaptivity. It is found that, although part of the domain adapts to a lower polynomial order through the adaptive algorithm, the calculated KE remains the same as the results obtained using the highest polynomial order.

\subsection{Compressible reacting flows}

The proposed method is finally applied to the simulation of compressible reacting flows. The compressible Navier-Stokes equations for a mixture of \(N_S\) miscible reacting chemical species are solved with the $p$-adaptive LABFM. The governing equations, written in index notation with Greek subscripts denoting spatial coordinates and Latin subscripts indicating physical entities such as phases or species for clarity, are given by,

\begin{align}
\frac{\partial \rho}{\partial t} + \frac{\partial \rho u_\alpha}{\partial x_\alpha} &= 0  \\
\frac{\partial \rho u_\alpha}{\partial t} + \frac{\partial \rho u_\alpha u_\beta}{\partial x_\beta} &= -\frac{\partial p_f}{\partial x_\alpha} + \frac{\partial \tau_{\alpha\beta}}{\partial x_\beta} + \rho f_\alpha \\
\frac{\partial \rho E}{\partial t} + \frac{\partial \rho u_\alpha E}{\partial x_\alpha} &= -\frac{\partial p_f u_\alpha}{\partial x_\alpha} - \frac{\partial q_\alpha}{\partial x_\alpha} + \frac{\partial \tau_{\alpha\beta} u_\alpha}{\partial x_\beta} + \rho \sum_{k=1}^{N_S} Y_k f_\alpha(u_\alpha - V_{\alpha.k}) \\
\frac{\partial \rho Y_k}{\partial t} + \frac{\partial \rho u_\alpha Y_k}{\partial x_\alpha} &= \dot\omega_k - \frac{\partial \rho V_{\alpha.k} Y_k}{\partial x_\alpha}, \quad \forall k \in [1, N_S] 
\end{align}
where $u_\alpha$ is the velocity component in the $\alpha$-direction, \(\tau_{\alpha\beta}\) is the component of viscous stress tensor, \(f_\alpha\) is the body force, \(q_\alpha\) is the heat flux vector, \(Y_k\) is the mass fraction of species \(k \in [1,N_S]\), \(V_{\alpha.k}\) is the molecular diffusion velocity of species \(k\) in the $\alpha$-direction, and \(\dot\omega_k \) is the production rate of species \(k\). The total energy \(E\) has a relationship with other thermodynamic quantities via the caloric equation of state according to,
\begin{equation}
E = \sum_{k} h_{k} Y_{k} - \frac{p_f}{\rho} + \frac{1}{2} u_{\alpha} u_{\alpha}
\end{equation}
where \(h_k\) is the enthalpy of species \(k\).

The governing equations are closed by the thermal equation of state
\begin{equation}
p_f = \rho R_0T\sum_{k} Y_k/W_k
\end{equation}
where \(R_0\) is the universal gas constant, \(T\) is the temperature, \(W_k\) is the molar mass of species \(k\).

The temperature dependencies of thermodynamic quantities (i.e., heat capacity \(c_{p,k}\), enthalpy \(h_k\)) for each species follow polynomial forms fitted to standard NASA polynomials \cite{gordon1976computer}. Two reaction mechanisms are utilised in this work for hydrogen flames, including the 21-step, 9-species hydrogen-air mechanism of~\cite{li2004updated}, and a single-step irreversible Arrhenius mechanism as in~\cite{dominguez2023stable}. The latter is constructed for dilute hydrogen combustion and is tuned to match the laminar flame speed of hydrogen at atmospheric pressure and an equivalence ratio, i.e., the actual ratio of fuel-to-oxidiser in a mixture to the stoichiometric ratio, of 0.32. For this mechanism, the temperature dependence of transport properties is prescribed via a power law with an exponent of 0.7, and a constant Lewis number approximation is applied for molecular diffusion, allowing the thermal conductivity and molecular diffusivities to be determined from the definitions of the Prandtl and Lewis numbers. For the multi-step reaction mechanism, a mixture-averaged model for transport properties is employed, following the combination rules of~\cite{ern1994multicomponent, ern1995fast}, with the Hirschfelder-Curtiss approximation for molecular diffusion~\cite{hirschfelder1964molecular}. Soret and Dufour effects are neglected, consistent with~\cite{howarth2022empirical}. Time integration is performed using an explicit third-order Runge-Kutta scheme. The implementation of $p$-adaptive LABFM for compressible reacting flow simulation, as described above, is achieved within the message passing interface (MPI)-enabled DNS combustion code, the SUNSET code (\textbf{S}calable \textbf{U}nstructured \textbf{N}ode-\textbf{SET} code)~\cite{sunset}. A detailed description of the mathematical model and numerical implementation is available in~\cite{king2024cmame}.

\subsubsection{Freely propagating $H_{2}$-air flames}

This case simulates planar, laminar, freely propagating stoichiometric hydrogen-air flames. A rectangular two-dimensional domain is considered, with length $L_{x}=L=10mm$ and a narrow aspect ratio, $L_{y}=L/40$. The domain is periodic in $y$, with inflow and outflow boundaries in $x$. The temperature at the inflow is set to $T_{in}=300K$, the velocity to $U=2.05m/s$, and mass fractions corresponding to a stoichiometric hydrogen-air mixture are imposed. A partially non-reflecting outflow condition \cite{king2024cmame} is used to maintain an outflow pressure of $p_{f.out}=10^{5}Pa$. The 9-species, 21-step reaction mechanism of~\cite{li2004updated} is employed.

\begin{figure}[h]
    \centering
    \includegraphics[width=0.8\textwidth]{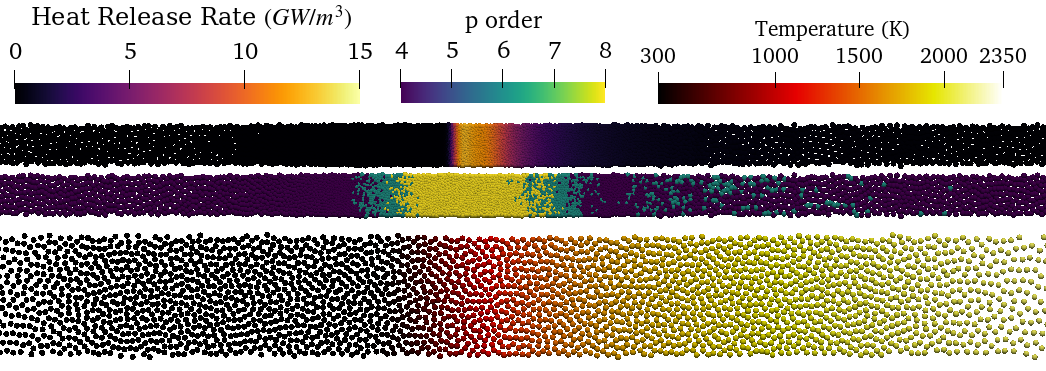}
    \caption{Results and node distribution for the planar laminar hydrogen-air flame simulations: heat release rate (top), polynomial order $p$ (middle), and temperature (bottom)}\label{fig: laminar_flame}
\end{figure}

\begin{figure}[h]
    \centering
    \includegraphics[width=0.49\textwidth]{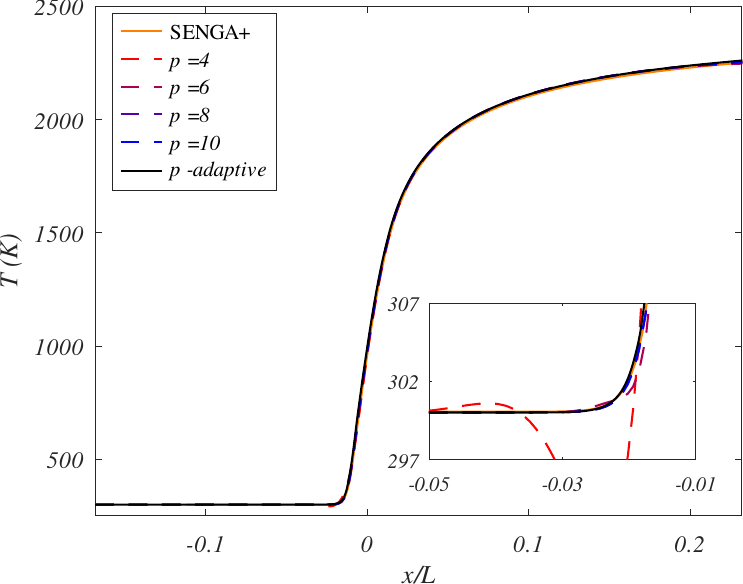}
    \includegraphics[width=0.49\textwidth]{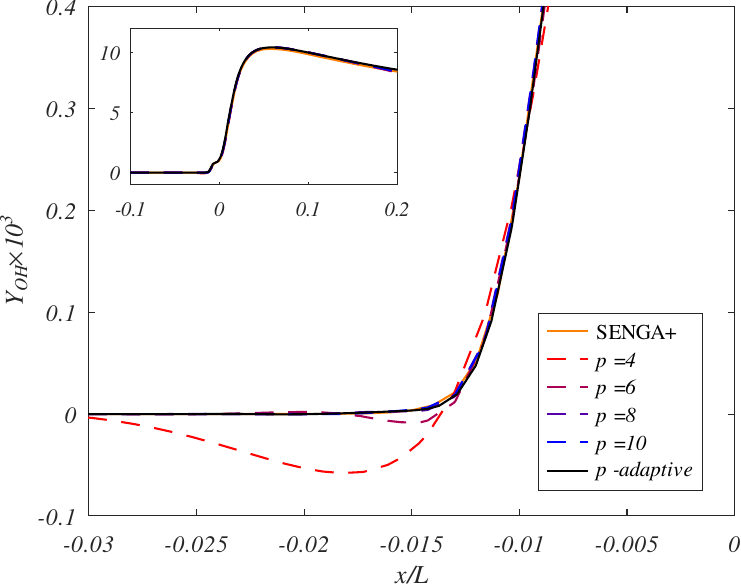}
    \caption{Spatial profiles of temperature (left) and $\mathrm{OH}$ mass fraction (right) for a freely propagating stoichiometric hydrogen-air flame, obtained with the present adaptive code (black lines), with the non-adaptive version at orders $p\in\{4,6,8,10\}$ (dashed lines), and SENGA+ (orange line). In all cases, the resolution in the flame region is $L/1000$.\label{fig:flamedisc}} 
\end{figure}

\begin{figure}[!htbp]
    \centering
    \includegraphics[width=0.99\textwidth]{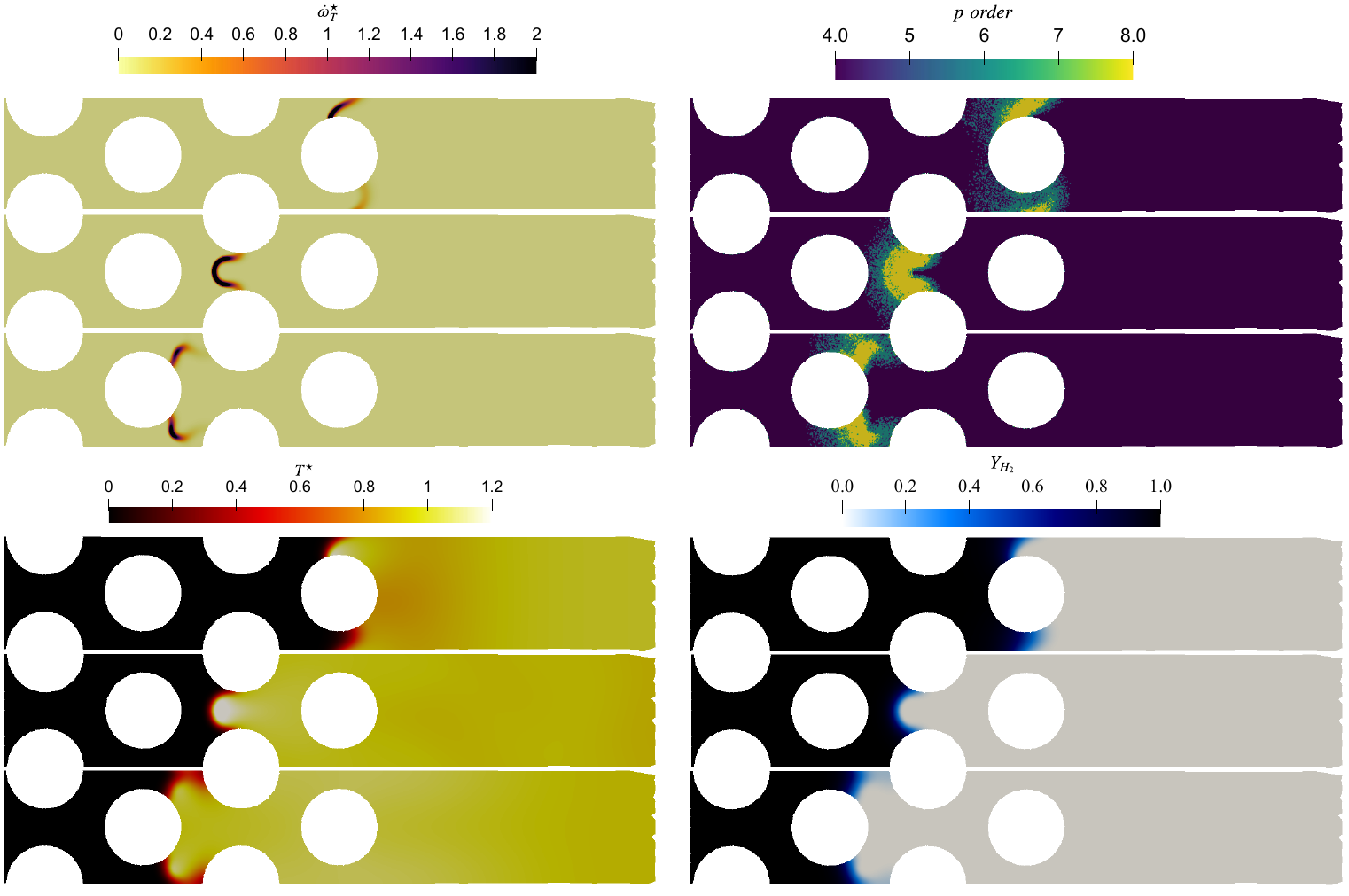}
    \caption{Time evolution of a hydrogen flame propagating through a porous matrix, showing the heat release rate $\dot{\omega}^{\star}_{T}$ (top left), the local polynomial order $p$ (top right), the temperature $T^*$ (bottom left), and the mass fraction $Y_{H_2}$ (bottom right). Within each panel, the first, second, and third rows correspond to $t^{\star}=11$, $t^{\star}=18$, and $t^{\star}=23$, respectively.\label{fig:adaptive1}} 
\end{figure}

For the discretisation, $s/L=0.003$ is set at the inflow and outflow boundaries, reduced linearly over a region of $0.1L$ to $s/L=0.001$ in the central $0.1L$ of the domain, corresponding to a resolution of $s=10\,\mu\mathrm{m}$. This configuration ensures that the flame structure is resolved by at least 10 nodes. $\phi$ in the refiement procedures is taken as the mass fraction of $Y_{H_2}$. The upper and lower thresholds for adaptivity are taken as \(\varepsilon_U=10^{-3}\) and \(\varepsilon_L=10^{-4}\). The simulation is run until a steady flame profile is achieved.

Figure~\ref{fig:flamedisc} shows the spatial profiles of temperature (left) and the mass fraction of the hydroxyl radical $\mathrm{OH}$ (right) for a freely propagating hydrogen-air flame, obtained using the present adaptive method, the non-adaptive version at orders $p\in\{4,6,8,10\}$, and the reference code SENGA+~\cite{cant_2012}. SENGA+ is one of the leading codes for combustion DNS based on a tenth-order finite difference scheme with a fourth-order explicit Runge–Kutta method for time integration. Localised differences are observed for $p = 4$ and $6$ compared to other results. A possible reason could be the lower resolving power at lower polynomial orders, resulting in insufficient resolution to accurately compute the advection terms. This leads to slight oscillatory behaviour just upstream of the flame front, with local minima in temperature and regions of unphysical negative $\mathrm{OH}$ mass fraction. For $p\geq8$, these oscillations are absent, and the results closely match the reference solution. Notably, the $p$-adaptive scheme, which varies the order $p$ from $4$ to $8$, eliminates the issues associated with uniformly low-order schemes while achieving accuracy comparable to the uniform $p=8$ case. In addition, the $p$-adaptive code achieves a speed-up of approximately 15–20\% compared to the uniform $p=8$ case. Note that the current implementation is based on an MPI-enabled parallel scheme with static domain decomposition, indicating promising potential for further increase of computational speed-up with the implementation of dynamic load balancing.

\subsubsection{Flame propagation through a cylinder array}

In this case, the adaptive scheme is applied to the analysis of a hydrogen flame propagating through a porous matrix constructed from an array of cylinders. This scenario features complex geometries and dynamic flame front evolution. The propagation of the flame is investigated within an idealised porous medium represented by a hexagonal array of cylinders. 

The hydrogen-air mixture under consideration has an equivalence ratio of \(\phi = 0.32\), at atmospheric pressure and a temperature of \(300\,\mathrm{K}\). The laminar flame speed and thermal flame thickness are \(S_L \approx 0.07\,\mathrm{m\,s^{-1}}\) and \(\delta_L \approx 1\,\mathrm{mm}\), respectively, serving as reference scales for non-dimensionalisation. The cylinder diameter is fixed at \(D/\delta_L = 5\), while the cylinder spacing is taken as \(S/D = 5/3 \) . The computational domain consists of a rectangular region with a length of $L=10D$ and a width of $S$. Lateral boundaries are periodic, while streamwise boundaries employ non-reflecting inflow and outflow conditions, as described in \cite{king2024cmame}. The temperature at the inflow is set to $T_{in}=300\,\mathrm{K}$, while the outflow boundary tracks a pressure of $p_{f,out}=10^{5}\,\mathrm{Pa}$. Adiabatic no-slip conditions are imposed on the cylinder surfaces. The domain is discretised with a non-uniform resolution, varying from \(s_{min}= D/80\) at the cylinder surface, to \(s_{out}= 6s_{min}\) at the outflow, and \(s_{in}= s_{min}\) at the inflow. As in the previous test, the upper and lower thresholds for adaptivity are set as \(\varepsilon_U=10^{-3}\) and \(\varepsilon_L=10^{-4}\). The simulations are initialised with a laminar flame positioned just downstream of the cylinder array and perturbed by a small, random, multi-mode disturbance with an amplitude of \(\delta_L/4\).

Fig.~\ref{fig:adaptive1} shows the contour plot of the heat release rate, the temperature, the mass fraction of $H_2$, and the local polynomial order during the flame propagation through pore structures at different time instants. The flame propagates relatively evenly through each pore,  with a uniform evolution of $Y_{H_2}$. An enhanced heat release rate is observed as the flame passes through the pore throat, alongside a local increase of temperature. The dynamic evolution of the polynomial order aligns well with the flame's propagation. This demonstrates that the adaptive algorithm effectively captures the flame front, a region where high-order accuracy is necessary, as it propagates through a complex geometry.

\begin{figure}[ht]
    \centering
    \includegraphics[width=0.49\textwidth]{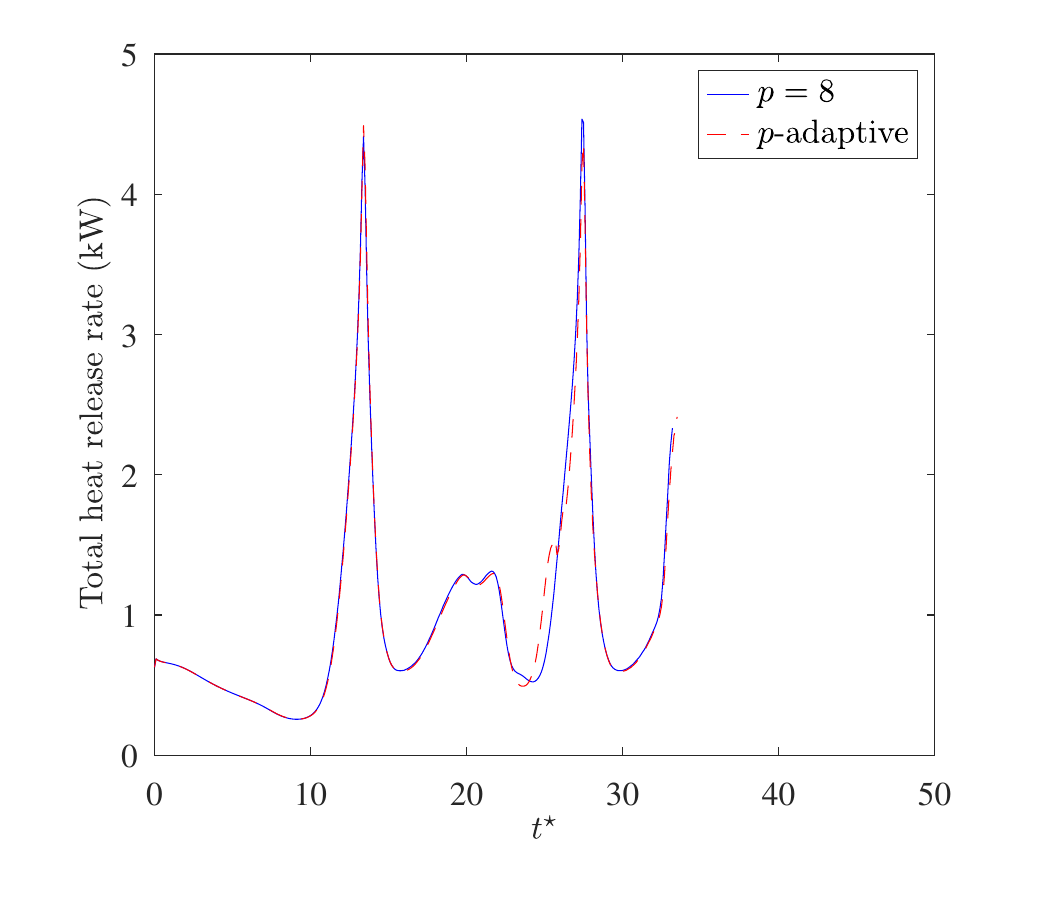}
    \includegraphics[width=0.49\textwidth]{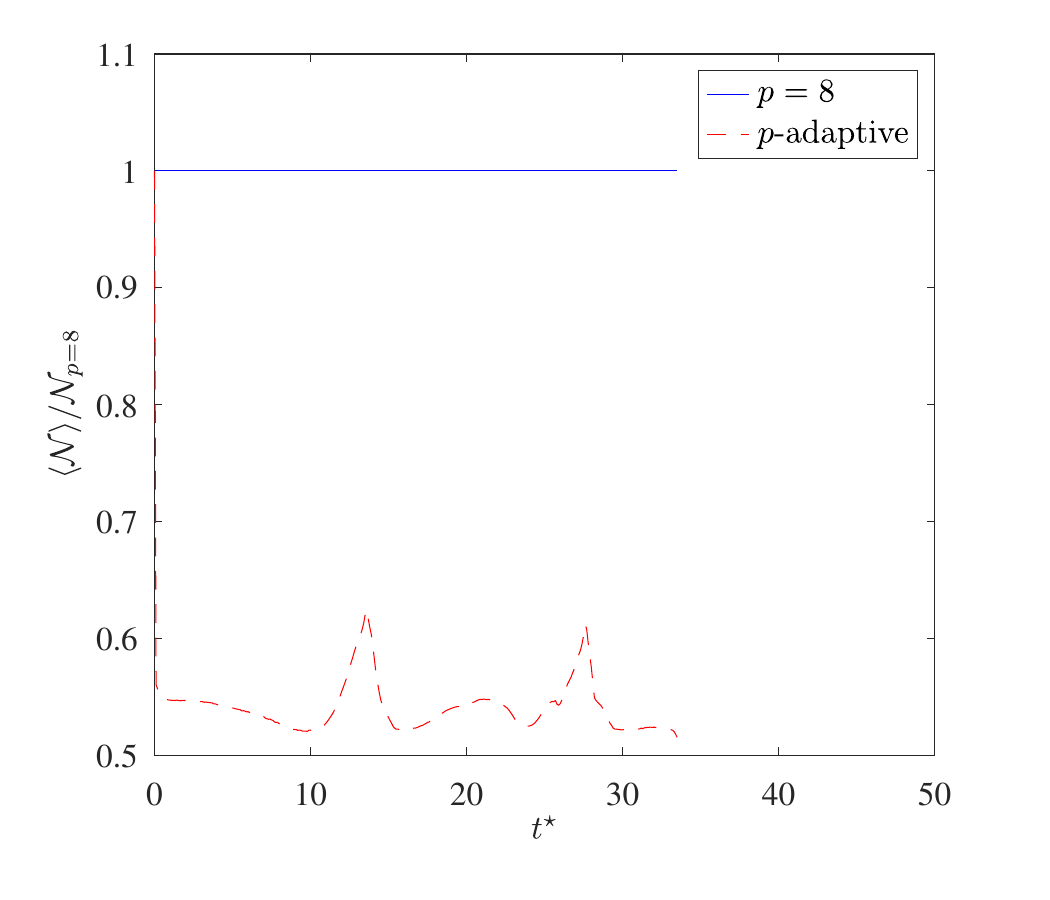}
    \caption{Comparison of results obtained using $p$-adaptivity and non-adaptivity with $p=8$, showing time evolution of the integrated heat release rate (left) and the normalised neighbour count (right).\label{fig:adaptive2}} 
\end{figure}

The time evolution of the integrated heat release rate and the normalised neighbour counts per node, comparing fixed $p=8$ and $p$-adaptivity, is shown in Fig.~\ref{fig:adaptive2}. The heat release fluctuates as the flame propagates through the pores, while the adaptive algorithm captures the varied heat release and adjusts accordingly. The results for the heat release rate obtained using the $p$-adaptive method are nearly identical to those obtained with $p=8$, although a small discrepancy occurs at around $t^* = 23$. This could be attributed to the thermal diffusivity being sensitive to discretisation error. However, the computational cost for mesh-free interpolation, represented by the neighbour counts per node, is reduced by nearly half with the adaptive algorithm. This demonstrates the improved computational efficiency and high accuracy possible with the proposed $p$-adaptive mesh-free approach.

\section{Conclusions}
\label{sec:conclusions}

In this work, a $p$-adaptive, high-order mesh-free scheme based on the Local Anisotropic Basis Function Method (LABFM) is developed for flow problems in complex geometries. A refinement indicator based on the mesh-free estimate of local discretisation error of the Laplacian operator is formulated, and a corresponding refinement strategy is investigated to increase/decrease the order of the LABFM polynomial reconstruction at local computational nodes. The performance of the new $p$-adaptive mesh-free scheme is tested through several benchmark tests and compared with the original LABFM with a uniform polynomial order. The results demonstrate that the developed $p$-adaptive LABFM effectively captures regions where high accuracy is necessary and adjusts the approximation order accordingly. This approach leads to improved solution stability compared to a fixed low-order method, while achieving comparable accuracy to a fixed high-order method with better computational efficiency. The new mesh-free $p$-adaptive algorithm is then applied to resolve flame-structure interactions in complex geometries, showing a notable reduction in computational cost compared to the non-adaptive high-order implementation while maintaining high accuracy.

However, it is also worth mentioning that when extending the proposed method to parallel computation, there could be an issue of unbalanced load distribution due to the adaptive algorithm. In such a situation, a dynamic domain partition or task-based parallelism may be needed to realise the potential speedup. A parallel $p$-adaptive LABFM with a dynamic load-balancing technique will be a subject of future investigation.

\begin{acknowledgments}
This work was funded by the Engineering and Physical Sciences Research Council (EPSRC) grant EP/W005247/1 and EP/W005247/2. JK is funded by the Royal Society via a University Research Fellowship (URF\textbackslash R1\textbackslash 221290). The authors would like to acknowledge the use of the Computational Shared Facility at The University of Manchester.
\end{acknowledgments}





\nocite{*}

\bibliography{p_adapt_LABFM}

\end{document}